\documentclass[11pt]{amsart}
\usepackage[T1]{fontenc}
\usepackage[utf8]{inputenc}
\usepackage{lmodern}
\usepackage{amsmath,amssymb,amsthm,mathtools}
\usepackage[a4paper,margin=1in]{geometry}
\usepackage{enumitem}
\usepackage{tikz}
\usepackage{eucal}
\usepackage{float}
\usepackage{placeins}
\usepackage{tabularx}
\usetikzlibrary{calc}
\usepackage{booktabs}
\usepackage[hidelinks]{hyperref}

\numberwithin{equation}{section}

\newtheorem{theorem}{Theorem}[section]
\newtheorem{lemma}[theorem]{Lemma}
\newtheorem{proposition}[theorem]{Proposition}

\theoremstyle{definition}
\newtheorem{definition}[theorem]{Definition}

\newtheorem{remark}[theorem]{Remark}

\newcommand{\nmax}{64}
\newcommand{\R}{\mathbb R}
\newcommand{\Acal}{\mathcal A}
\newcommand{\Kcal}{\mathcal K}
\newcommand{\Ical}{\mathcal I}
\DeclareMathOperator{\rank}{rank}

\title{Area Stability for Largest Small Polygons}
\author{Jos\'e Ayala Hoffmann}
\address{Universidad de Tarapac\'a, Iquique, Chile}
\email{jayalhoff@gmail.com}

\author{Joaqu\'in L\'opez Guti\'errez}
\address{Universidad de Tarapac\'a, Iquique, Chile}
\email{juakinigna@gmail.com}
\date{}



\begin{document}

\begin{abstract} \baselineskip=20 true pt
\baselineskip=1.1\normalbaselineskip
We develop a second order calculus for area maximisation among
small polygons of diameter at most one, culminating in rigorous
computational proofs of stationarity, spectral stability and
strict local maximality modulo orientation preserving Euclidean
motions for the polygons \(Q_n\), whose diameter graphs consist
of a cycle on \(n-1\) vertices and one pendant edge,
for every even \(n\) with \(6\le n\le64\).

At a stationary configuration, positive multipliers force
first order area decrease along feasible directions that
shorten at least one active diameter. Along directions that preserve all active diameters to first
order, the reduced Hessian governs second order
area change modulo Euclidean motions. Positive multipliers,
independent active constraint gradients, and a negative reduced
spectrum imply strict local maximality. Outward rounded interval arithmetic establishes existence of
strictly convex stationary configurations \(Q_n\) with all \(n-3\) reduced eigenvalues strictly negative. Together, these results yield explicit neighbourhood radii
and rigorous quadratic area decrease bounds for all
sufficiently nearby feasible polygons. 

The reduced spectrum also detects second order degeneracy
and instability, as illustrated by an optimal quadrilateral
and the regular hexagon.
\end{abstract}
	\maketitle

    \baselineskip=20 true pt
\baselineskip=1.1\normalbaselineskip
	\section{Prelude}
What is the largest possible area of a convex $n$-gon whose
vertices are pairwise at distance at most one?
Such polygons are called small polygons, and the corresponding
maximisation problem is known as the largest area small polygon problem. 

Reinhardt in 1922 proved that the regular small
	\(n\)-gon is optimal when \(n\) is odd, whereas it is not optimal for
	even \(n\ge6\) \cite{reinhardt}. Graham determined the largest small
	hexagon \cite{graham}, and Audet et al. resolved the octagon case
	\cite{audethansenmessinexiong}. Foster and Szab\'o then established a 
    structural result valid for every even \(n\ge6\). The diameter graph of an area maximiser, whose edges join pairs of vertices at distance \(1\), consists of a
    cycle on \(n-1\) vertices with one pendant edge \cite{fosterszabo}. We call
    this a Graham type diameter graph. 
	
Mossinghoff developed explicit constructions and area bounds
\cite{mossinghoff}, subsequently improved by Bingane
\cite{binganedcg} and by Bingane and Mossinghoff
\cite{binganemossinghoff}.
Numerical optimisation has also provided candidate polygons,
with Henrion and Messine computing candidates for the largest
small decagon and dodecagon by global polynomial optimisation
and reporting numerical evidence of global optimality
\cite{henrionmessine}, while Pintér, Kampas and
Castillo~\cite{pinterkampascastillo} obtained further numerical
candidates using nonlinear programming. These identify promising configurations, but numerical evidence alone does not guarantee optimality.
    Once a stationary candidate has been found, we may ask how its area changes under perturbations.
    
	A diameter constraint is active when the corresponding pair of vertices is at distance one. At a stationary configuration, feasible first order directions
either shorten at least one active diameter or preserve all
active diameters to first order. If a direction shortens at
    least one of them, positive active multipliers imply a strictly negative first 
    variation of area. The directions preserving all active diameters to first order
form the active equality space. On this space, the first variation
vanishes, and we examine the Hessian of the Lagrangian. Restricting this Hessian to the active equality space orthogonal
to infinitesimal Euclidean motions gives the reduced Hessian,
whose eigenvalues form the reduced spectrum. When the
    active constraint gradients are independent and all active multipliers are 
    positive, negative definiteness of the reduced Hessian gives a second order
    sufficient condition for strict local maximality modulo orientation preserving
    Euclidean motions \cite{nocedalwright}. We study this spectrum and turn the 
    criterion into explicit bounds on area loss.
	
The quadratic structure of the problem allows us to make this criterion
quantitative. Both the area \(\Acal\) and the constraint functions
\(g_{ij}(x)=\lVert x_i-x_j\rVert^2-1\) are quadratic, so the Lagrangian
expansion at a stationary configuration \(c\) is exact. 
Under the rank, multiplier and spectral conditions above,
Theorem~\ref{thm:quant-stability} gives explicit constants \(\rho,\kappa>0\) such that every configuration \(x\) satisfying
\(\lVert x_i-x_j\rVert\leq1\) for all
\(\{i,j\}\in E(c)\) is such that
$
\Acal(x)\leq\Acal(c)-\kappa\operatorname{dist}(x,c)^2$
whenever $\operatorname{dist}(x,c)\leq\rho$. Here \(E(c)\) is the set of active diameter pairs at \(c\), and
\(\operatorname{dist}(x,c)\) is the distance from \(x\) to the orbit
of \(c\) under orientation preserving Euclidean motions. The
    constants depend explicitly on bounds for the smallest row singular value of
    the active constraint matrix, the active multipliers and the largest reduced 
    eigenvalue, together with the maximum degree of the active diameter graph. The
    estimate also covers configurations that shorten active diameters and requires
    no inactive distance margin.

    We apply our criterion to stationary configurations with a
    Graham type diameter graph. Following the notation
    of Mossinghoff~\cite{mossinghoff}, we write \(Q_n\) for the
    configuration defined in Theorem~\ref{thm:main}. For each even
    \(n\) with \(6\le n\le\nmax\), outward rounded interval arithmetic
    establishes the existence of an exact, strictly convex stationary
    configuration \(Q_n\), with independent active constraint gradients
    and strictly positive multipliers. The entire reduced spectrum of \(Q_n\) is enclosed in 
    \(n-3\) pairwise disjoint negative intervals, each containing one eigenvalue of multiplicity one. Since the largest reduced eigenvalue has multiplicity one,
the unit direction of weakest quadratic area decrease in the
reduced tangent space is unique up to sign. This eigenvalue also enters the explicit stability constants
    \(\rho_n, \kappa_n\). The quantitative estimate
    then establishes that \(Q_n\) is a strict local maximiser modulo
    orientation preserving Euclidean motions, with a certified bound
    on area loss for all nearby feasible configurations. The computations
    also enclose \(\Acal(Q_n)\) in an interval of width less than
    \(10^{-11}\).

    Beyond the strictly stable configurations \(Q_n\), the
second order analysis also identifies degeneracy and instability.
The quadrilateral \(Q_4\) belongs to a nonisolated family of area
maximisers. On the equality space defined by its positive
multipliers, removal of rigid motions leaves two zero eigenvalues
corresponding to the parameters of this family. The remaining
active diameter inequalities restrict its feasible directions.
For the regular hexagon, two positive reduced eigenvalues
violate the second order necessary condition for local maximality.
	
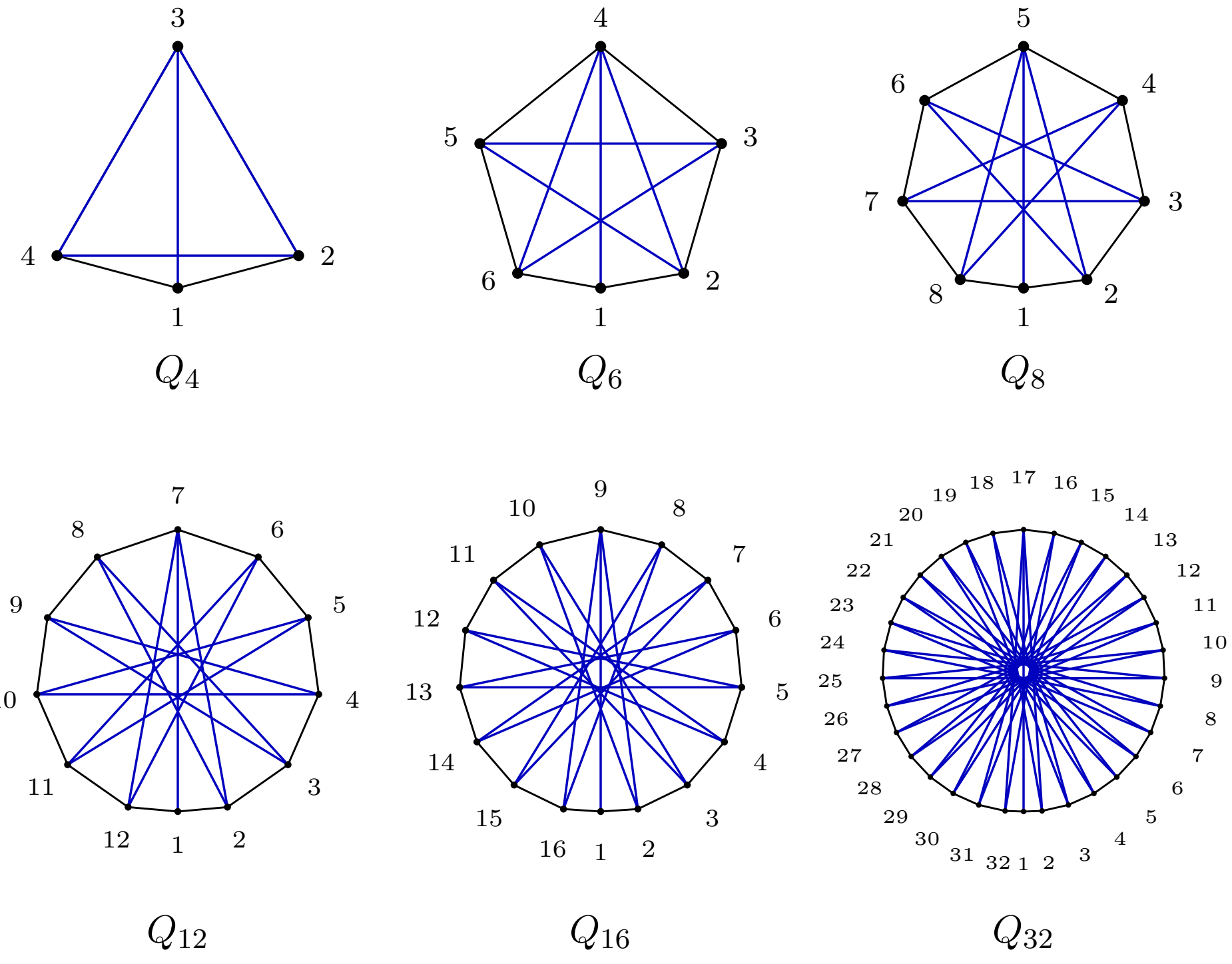
\begin{figure}[!ht]
		\centering
		\begin{tikzpicture}[
			scale=1.00,
			vertex/.style={circle,fill=black,inner sep=1.1pt},
			edge/.style={blue!75!black,line width=0.65pt},
			cycleedge/.style={black,line width=0.55pt}
			]
			
			
			\begin{scope}[shift={(0,0)},scale=2.4]
				\coordinate (A1) at (0,0);
				\coordinate (A2) at (0.500000,0.133975);
				\coordinate (A3) at (0,1.000000);
				\coordinate (A4) at (-0.500000,0.133975);
				\draw[cycleedge] (A1)--(A2)--(A3)--(A4)--cycle;
				\foreach \i/\j in {1/3,2/3,2/4,3/4}{\draw[edge] (A\i)--(A\j);}
				\foreach \i/\dx/\dy in {1/0.00/-0.12,2/0.12/0.00,3/0.00/0.12,4/-0.12/0.00}{
					\node[vertex] at (A\i) {};
					\node at ($(A\i)+(\dx,\dy)$) {\scriptsize \(\i\)};
				}
				\node at (0,-0.35) {\(Q_4\)};
			\end{scope}
			
			\begin{scope}[shift={(4.2,0)},scale=2.4]
				\coordinate (B1) at (0,0);
				\coordinate (B2) at (0.343771,0.060947);
				\coordinate (B3) at (0.500000,0.597649);
				\coordinate (B4) at (0,1.000000);
				\coordinate (B5) at (-0.500000,0.597649);
				\coordinate (B6) at (-0.343771,0.060947);
				\draw[cycleedge] (B1)--(B2)--(B3)--(B4)--(B5)--(B6)--cycle;
				\foreach \i/\j in {1/4,2/4,2/5,3/5,3/6,4/6}{\draw[edge] (B\i)--(B\j);}
				\foreach \i/\dx/\dy in {1/0.00/-0.12,2/0.12/-0.03,3/0.12/0.03,4/0.00/0.12,5/-0.12/0.03,6/-0.12/-0.03}{
					\node[vertex] at (B\i) {};
					\node at ($(B\i)+(\dx,\dy)$) {\scriptsize \(\i\)};
				}
				\node at (0,-0.35) {\(Q_6\)};
			\end{scope}
			
			\begin{scope}[shift={(8.4,0)},scale=2.4]
				\coordinate (C1) at (0,0);
				\coordinate (C2) at (0.262142,0.034971);
				\coordinate (C3) at (0.500000,0.359621);
				\coordinate (C4) at (0.409092,0.776216);
				\coordinate (C5) at (0,1.000000);
				\coordinate (C6) at (-0.409092,0.776216);
				\coordinate (C7) at (-0.500000,0.359621);
				\coordinate (C8) at (-0.262142,0.034971);
				\draw[cycleedge] (C1)--(C2)--(C3)--(C4)--(C5)--(C6)--(C7)--(C8)--cycle;
				\foreach \i/\j in {1/5,2/5,2/6,3/6,3/7,4/7,4/8,5/8}{\draw[edge] (C\i)--(C\j);}
				\foreach \i/\dx/\dy in {1/0.00/-0.12,2/0.10/-0.06,3/0.13/0.00,4/0.11/0.07,5/0.00/0.12,6/-0.11/0.07,7/-0.13/0.00,8/-0.10/-0.06}{
					\node[vertex] at (C\i) {};
					\node at ($(C\i)+(\dx,\dy)$) {\scriptsize \(\i\)};
				}
				\node at (0,-0.35) {\(Q_8\)};
			\end{scope}


            \begin{scope}[shift={(0,-5.2)},scale=2.8]
                \coordinate (D1) at (0.000000000,0.000000000);
                \coordinate (D2) at (0.176161368,0.015638698);
                \coordinate (D3) at (0.390889812,0.166088565);
                \coordinate (D4) at (0.500000000,0.415725636);
                \coordinate (D5) at (0.462314079,0.687666096);
                \coordinate (D6) at (0.285342781,0.902776749);
                \coordinate (D7) at (0.000000000,1.000000000);
                \coordinate (D8) at (-0.285342781,0.902776749);
                \coordinate (D9) at (-0.462314079,0.687666096);
                \coordinate (D10) at (-0.500000000,0.415725636);
                \coordinate (D11) at (-0.390889812,0.166088565);
                \coordinate (D12) at (-0.176161368,0.015638698);
                \foreach \i in {1,...,11}{
                    \pgfmathtruncatemacro{\j}{\i+1}
                    \draw[cycleedge] (D\i)--(D\j);
                }
                \draw[cycleedge] (D12)--(D1);
                \foreach \i/\j in {1/7,2/7,2/8,3/8,3/9,4/9,4/10,5/10,5/11,6/11,6/12,7/12}{\draw[edge] (D\i)--(D\j);}
                \coordinate (DO) at (0,0.5);
                \foreach \i in {1,...,12}{
                    \node[vertex,inner sep=0.7pt] at (D\i) {};
                    \node[font=\tiny,inner sep=0pt] at ($(DO)!1.24!(D\i)$) {\i};
                }
                \node at (0,-0.43) {$Q_{12}$};
            \end{scope}

            \begin{scope}[shift={(4.2,-5.2)},scale=2.8]
                \coordinate (E1) at (0.000000000,0.000000000);
                \coordinate (E2) at (0.132041376,0.008755794);
                \coordinate (E3) at (0.307035216,0.094013608);
                \coordinate (E4) at (0.438987437,0.246881396);
                \coordinate (E5) at (0.500000000,0.440505245);
                \coordinate (E6) at (0.479349059,0.642681951);
                \coordinate (E7) at (0.380114368,0.820529574);
                \coordinate (E8) at (0.216369963,0.946097535);
                \coordinate (E9) at (0.000000000,1.000000000);
                \coordinate (E10) at (-0.216369963,0.946097535);
                \coordinate (E11) at (-0.380114368,0.820529574);
                \coordinate (E12) at (-0.479349059,0.642681951);
                \coordinate (E13) at (-0.500000000,0.440505245);
                \coordinate (E14) at (-0.438987437,0.246881396);
                \coordinate (E15) at (-0.307035216,0.094013608);
                \coordinate (E16) at (-0.132041376,0.008755794);
                \foreach \i in {1,...,15}{
                    \pgfmathtruncatemacro{\j}{\i+1}
                    \draw[cycleedge] (E\i)--(E\j);
                }
                \draw[cycleedge] (E16)--(E1);
                \foreach \i/\j in {1/9,2/9,2/10,3/10,3/11,4/11,4/12,5/12,5/13,6/13,6/14,7/14,7/15,8/15,8/16,9/16}{\draw[edge] (E\i)--(E\j);}
                \coordinate (EO) at (0,0.5);
                \foreach \i in {1,...,16}{
                    \node[vertex,inner sep=0.7pt] at (E\i) {};
                    \node[font=\tiny,inner sep=0pt] at ($(EO)!1.29!(E\i)$) {\i};
                }
                \node at (0,-0.43) {$Q_{16}$};
            \end{scope}

            \begin{scope}[shift={(8.4,-5.2)},scale=2.8]
                \coordinate (F1) at (0.000000000,0.000000000);
                \coordinate (F2) at (0.065568370,0.002151921);
                \coordinate (F3) at (0.159015004,0.023341763);
                \coordinate (F4) at (0.249593069,0.064006612);
                \coordinate (F5) at (0.330705868,0.122120271);
                \coordinate (F6) at (0.398720235,0.195229587);
                \coordinate (F7) at (0.450875593,0.280393345);
                \coordinate (F8) at (0.485089623,0.374216363);
                \coordinate (F9) at (0.500000000,0.472963948);
                \coordinate (F10) at (0.495013293,0.572706350);
                \coordinate (F11) at (0.470327912,0.669474446);
                \coordinate (F12) at (0.426925924,0.759417870);
                \coordinate (F13) at (0.366532106,0.838960007);
                \coordinate (F14) at (0.291530427,0.904949750);
                \coordinate (F15) at (0.204757449,0.954829606);
                \coordinate (F16) at (0.108598809,0.986868019);
                \coordinate (F17) at (0.000000000,1.000000000);
                \coordinate (F18) at (-0.108598809,0.986868019);
                \coordinate (F19) at (-0.204757449,0.954829606);
                \coordinate (F20) at (-0.291530427,0.904949750);
                \coordinate (F21) at (-0.366532106,0.838960007);
                \coordinate (F22) at (-0.426925924,0.759417870);
                \coordinate (F23) at (-0.470327912,0.669474446);
                \coordinate (F24) at (-0.495013293,0.572706350);
                \coordinate (F25) at (-0.500000000,0.472963948);
                \coordinate (F26) at (-0.485089623,0.374216363);
                \coordinate (F27) at (-0.450875593,0.280393345);
                \coordinate (F28) at (-0.398720235,0.195229587);
                \coordinate (F29) at (-0.330705868,0.122120271);
                \coordinate (F30) at (-0.249593069,0.064006612);
                \coordinate (F31) at (-0.159015004,0.023341763);
                \coordinate (F32) at (-0.065568370,0.002151921);
                \foreach \i in {1,...,31}{
                    \pgfmathtruncatemacro{\j}{\i+1}
                    \draw[cycleedge] (F\i)--(F\j);
                }
                \draw[cycleedge] (F32)--(F1);
                \foreach \i/\j in {1/17,2/17,2/18,3/18,3/19,4/19,4/20,5/20,5/21,6/21,6/22,7/22,7/23,8/23,8/24,9/24,9/25,10/25,10/26,11/26,11/27,12/27,12/28,13/28,13/29,14/29,14/30,15/30,15/31,16/31,16/32,17/32}{\draw[edge] (F\i)--(F\j);}
                \coordinate (FO) at (0,0.5);
                \foreach \i in {1,...,32}{
                    \node[vertex,inner sep=0.55pt] at (F\i) {};
                    \node[font=\fontsize{5.4}{6}\selectfont,inner sep=0pt] at ($(FO)!1.37!(F\i)$) {\i};
                }
                \node at (0,-0.43) {$Q_{32}$};
            \end{scope}

        \end{tikzpicture}
        \caption{Examples of locally area maximising small polygons. Black segments show boundary edges outside the declared edge set;
        blue segments mark the declared unit distance constraints.
        Vertex numbers follow the manuscript indexing.      The drawing of $Q_4$ uses the configuration of
        Proposition~\ref{prop:q4}. The other panels use nominal
        coordinates from the LSPectra inputs.}
        \label{fig:diameter-graphs-all}
    \end{figure}

The interval computations are carried out by our Julia program LSPectra~\cite{lspectra}. Figure~\ref{fig:diameter-graphs-all} shows selected configurations
\(Q_n\) and their diameter graphs, with the vertex labelling used throughout the paper. Figure~\ref{fig:q64-diameter-framework}
shows the diameter graph of \(Q_{64}\), with its pendant edge highlighted.
	
	Sections~2--4 establish the first and second order theory
and derive quantitative stability bounds.
Section~\ref{sec:models} analyses degeneracy in the quadrilateral case.
Section~\ref{sec:lspectra} develops the computational verification
procedure, and Section~\ref{sec:main} establishes the main results
and presents the certified spectral and stability bounds.

	\section{Small polygons and diameter graphs}
	Throughout, \(\langle\cdot,\cdot\rangle\) denotes the standard Euclidean inner
	product on \(\R^2\) extended componentwise to \((\R^2)^n\).
	
	\begin{definition}
		\label{def:small-polygon}
		For \(n\geq 3\), a small \(n\)-gon is a labelled cyclic tuple
		\[
		c=(c_1,\ldots,c_n)\in(\R^2)^n
		\]
		satisfying
		\[
		\lVert c_i-c_j\rVert\leq 1
		\qquad
		\text{for all }i\neq j.
		\]
		Its area is
		\[
		\Acal(c)
		=
		\frac12\sum_{i=1}^n\det(c_i,c_{i+1}),
		\qquad
		c_{n+1}=c_1.
		\]
		It is positively oriented if \(\Acal(c)>0\).
	\end{definition}
	
\begin{remark}
\label{rem:shoelace}
The functional \(\Acal\) is signed area. For simple, positively oriented
polygons it corresponds to the enclosed area.
\end{remark}
	
	\begin{definition}\label{def:active-diameter-graph}
		Let \(c\) be a small \(n\)-gon. Its active diameter graph is
		\[
		G(c)=(\{1,\ldots,n\},E(c)),
		\]
		where
		\[
		E(c)=\{\{i,j\}:\lVert c_i-c_j\rVert =1\}.
		\]
		The edges of \(G(c)\) are called active diameter edges.
	\end{definition}
	
	\begin{definition}\label{def:fixed-active-diameter-class}
		Let \(c\) be a small \(n\)-gon. The fixed active diameter class at \(c\) is the
		subset of a neighbourhood of \(c\) in \((\R^2)^n\) cut out by
		\[
		\lVert x_i-x_j \rVert =1
		\qquad
		(\{i,j\}\in E(c)),
		\]
		\[
		\lVert x_i-x_j\rVert<1
		\qquad
		(\{i,j\}\notin E(c)),
		\]
		with the labelled cyclic order of the area formula held fixed.
	\end{definition}
	
	\begin{definition}\label{def:first-order-cone}
		The linearised feasible cone for the diameter inequalities at \(c\) is
		\[
		\Kcal(c)=
		\left\{
		v\in(\R^2)^n:
		\langle c_j-c_i,v_j-v_i\rangle\leq 0
		\text{ for every }\{i,j\}\in E(c)
		\right\}.
		\]
	\end{definition}
	
	\begin{definition}\label{def:equality-space}
		The active equality space is
		\[
		\Ical(c)=
		\left\{
		v\in(\R^2)^n:
		\langle c_j-c_i,v_j-v_i\rangle=0
		\text{ for every }\{i,j\}\in E(c)
		\right\}.
		\]
	\end{definition}
	
		Throughout, Euclidean motions mean orientation preserving Euclidean motions,
		acting without relabelling. Reflections reverse the signed area. Let
		\[
		J(x,y)=(-y,x).
		\]
		The infinitesimal rigid motion space at \(c\) is
		\[
		R(c)=
		\{(a+\omega Jc_1,\ldots,a+\omega Jc_n):a\in\R^2,\ \omega\in\R\}.
		\]
		Its orthogonal complement \(R(c)^\perp\) is taken with respect to the standard
		Euclidean inner product on \((\R^2)^n\).
	
	\begin{definition}
		Let \(c\) be a small \(n\)-gon. Its inactive distance margin is
		\[
		\delta(c)
		=
		\min_{\{i,j\}\notin E(c)}
		\bigl(1-\lVert c_i-c_j \rVert \bigr).
		\]
		The minimum is over unordered pairs of distinct indices and is defined as
		\(+\infty\) if there are no inactive pairs. Since the active graph contains
		all pairs at distance \(1\), this margin is positive for every small \(n\)-gon.
	\end{definition}

		For every small \(n\)-gon \(c\), \(R(c)\subseteq\Ical(c)\). Indeed, if
		\(v_i=a+\omega Jc_i\), then
		\[
		v_j-v_i=\omega J(c_j-c_i),
		\]
		and therefore
		\[
		\langle c_j-c_i,v_j-v_i\rangle
		=
		\omega\langle c_j-c_i,J(c_j-c_i)\rangle
		=
		0.
		\]

	\begin{remark}\label{rem:full-row-rank}
		For each active edge \(\{i,j\}\in E(c)\), set
		\[
		g_{ij}(x)=\lVert x_i-x_j\rVert^2-1.
		\]
		Order the rows of the active constraint matrix
		\[
		Dg(c)=\bigl(Dg_{ij}(c)\bigr)_{\{i,j\}\in E(c)}
		\]
		by a fixed enumeration of \(E(c)\). The same fixed rows define \(Dg(x)\)
		for nearby \(x\), regardless of its actual active set. In particular,
		\(\Ical(c)=\ker Dg(c)\). If \(Dg(c)\) has full row rank, then so does
		\(Dg(x)\) for every \(x\) in a neighbourhood of
		\(c\). The active equality set
		\[
		M=\{x:g_{ij}(x)=0\text{ for all }\{i,j\}\in E(c)\}
		\]
		is therefore a \(C^2\) submanifold near \(c\), with
		\[
		T_cM=\Ical(c).
		\]
		This is the regularity hypothesis required by the transverse slice maximality, Proposition~\ref{prop:slice-local-max}. We write
		\(\rank Dg(c)=|E(c)|\) when the full row rank condition holds.
	\end{remark}
	
	\begin{remark}\label{rem:contact-rate}
		For an active edge \(\{i,j\}\in E(c)\) and a variation
		\(v\in(\R^2)^n\),
		\[
		Dg_{ij}(c)[v]
		=2\langle c_j-c_i,v_j-v_i\rangle.
		\]
		Since \(\lVert c_j-c_i \rVert =1\),
		\[
		\frac{d}{dt}\bigg|_{t=0}
		\bigl\lVert(c_j+tv_j)-(c_i+tv_i)\bigr\rVert
		=\langle c_j-c_i,v_j-v_i\rangle
		=\frac12Dg_{ij}(c)[v].
		\]
		We call this scalar the signed diameter-contact rate. A negative rate moves the
		pair into the strict interior of the diameter constraint, a zero rate preserves
		the active contact to first order, and a positive rate points outside the
		feasible region. Consequently,
		\[
		\Kcal(c)=\{v:Dg(c)[v]\leq0\},
		\qquad
		\Ical(c)=\ker Dg(c),
		\]
		where the vector inequality is understood componentwise. The factor \(2\) is
		fixed by the unscaled convention \(g_{ij}(x)=\lVert x_i-x_j\rVert ^2-1\) and is retained
		throughout the multiplier and constrained Hessian formulae.
	\end{remark}
	\section{First variation and no ascent criteria}
	\begin{theorem}
		\label{thm:first-variation}
		For every \(c,v\in(\R^2)^n\) with fixed cyclic order,
		\[
		D\Acal(c)[v]
		=
		\frac12\sum_{i=1}^n\det(v_i,c_{i+1}-c_{i-1}),
		\]
		where indices are read cyclically. Equivalently,
		\[
		\nabla_{c_i}\Acal(c)
		=
		-\frac12J(c_{i+1}-c_{i-1}).
		\]
	\end{theorem}
	
	\begin{proof}
		The shoelace formula expresses the area as the polynomial
		\[
		\Acal(c)=\frac12\sum_{i=1}^n\det(c_i,c_{i+1}),
		\]
		and differentiating it along \(c+tv\) gives
		\[
		D\Acal(c)[v]
		=
		\frac12\sum_{i=1}^n
		\left[
		\det(v_i,c_{i+1})+\det(c_i,v_{i+1})
		\right].
		\]
		Reindexing the second sum,
		\[
		\sum_{i=1}^n\det(c_i,v_{i+1})
		=
		\sum_{i=1}^n\det(c_{i-1},v_i)
		=
		-\sum_{i=1}^n\det(v_i,c_{i-1}),
		\]
		so that
		\[
		D\Acal(c)[v]
		=
		\frac12\sum_{i=1}^n\det(v_i,c_{i+1}-c_{i-1}).
		\]
		The gradient formula follows from the identity \(\det(a,b)=\langle -Jb,a\rangle\).
	\end{proof}
	
	\begin{proposition}\label{prop:first-order-no-ascent}
		Let \(c\) be a positively oriented small \(n\)-gon with active diameter edge set
		\(E(c)\). Suppose there exist multipliers
		\[
		\lambda_{ij}\geq 0
		\qquad
		(\{i,j\}\in E(c)),
		\]
		with \(\lambda_{ij}=\lambda_{ji}\), such that
		\[
		\nabla_{c_i}\Acal(c)
		=
		2\sum_{j:\{i,j\}\in E(c)}
		\lambda_{ij}(c_i-c_j)
		\qquad
		\text{for every }i.
		\]
		Then \(D\Acal(c)[v]\leq 0\) for every \(v\in\Kcal(c)\). Thus \(c\) admits
		no first order area increasing direction on the linearised feasible cone.
		In particular, the first variation vanishes on the active equality space.
	\end{proposition}
	
	\begin{proof}
		The multiplier identity gives
		\[
		D\Acal(c)[v]
		=
		\sum_{i=1}^n\langle\nabla_{c_i}\Acal(c),v_i\rangle
		=
		2\sum_{i=1}^n
		\left\langle
		\sum_{j:\{i,j\}\in E(c)}
		\lambda_{ij}(c_i-c_j),
		v_i
		\right\rangle.
		\]
		Collecting the two vertex contributions from each unordered active edge yields
		\[
		D\Acal(c)[v]
		=
		2\sum_{\{i,j\}\in E(c)}
		\lambda_{ij}
		\left(
		\langle c_i-c_j,v_i\rangle+
		\langle c_j-c_i,v_j\rangle
		\right),
		\]
		that is,
		\[
		D\Acal(c)[v]
		=
		2\sum_{\{i,j\}\in E(c)}
		\lambda_{ij}
		\langle c_i-c_j,v_i-v_j\rangle.
		\]
		Since
		\[
		\langle c_i-c_j,v_i-v_j\rangle
		=
		\langle c_j-c_i,v_j-v_i\rangle
		\]
		and \(v\in\Kcal(c)\), every inner product is nonpositive. Since every
		\(\lambda_{ij}\geq 0\), every summand is nonpositive. Hence
		\(D\Acal(c)[v]\leq 0\).
	\end{proof}
	
	\begin{definition}
		A small \(n\)-gon \(c\) is stationary if the multiplier identity of
		Proposition~\ref{prop:first-order-no-ascent} holds for some
		\(\lambda\in\R^{E(c)}\), without a sign condition on its entries.
	\end{definition}
	
	\begin{lemma}\label{lem:positive-multipliers}
		Suppose the multipliers of Proposition~\ref{prop:first-order-no-ascent} are
		strictly positive,
		\[
		\lambda_{ij}>0
		\qquad
		(\{i,j\}\in E(c)).
		\]
		If \(v\in\Kcal(c)\) and \(D\Acal(c)[v]=0\), then
		\[
		\langle c_j-c_i,v_j-v_i\rangle=0
		\qquad
		(\{i,j\}\in E(c)).
		\]
		Thus every first order neutral direction in \(\Kcal(c)\) lies in
		\(\Ical(c)\).
	\end{lemma}
	
	\begin{proof}
		The proof of Proposition~\ref{prop:first-order-no-ascent} gives
		\[
		D\Acal(c)[v]
		=
		2\sum_{\{i,j\}\in E(c)}
		\lambda_{ij}
		\langle c_j-c_i,v_j-v_i\rangle.
		\]
		Each inner product is nonpositive because \(v\in\Kcal(c)\), and each
		multiplier is positive. A vanishing sum therefore forces every summand to
		vanish, so that
		\[
		\langle c_j-c_i,v_j-v_i\rangle=0
		\qquad
		(\{i,j\}\in E(c)).
		\]
	\end{proof}
	
The previous result identifies the directions that require a second order
analysis. With strictly positive multipliers, the first variation
vanishes on \(\Kcal(c)\) precisely on \(\Ical(c)=\ker Dg(c)\).
After removing infinitesimal rigid motions, we therefore examine
the second variation on \(\Ical(c)\cap R(c)^\perp\).

	\section{Second variation and local maximality}
	
	\begin{theorem}
		\label{thm:second-variation}
		For every \(c,v,w\in(\R^2)^n\), with indices taken modulo \(n\),
		\[
		D^2\Acal(c)[v,w]
		=
		\frac12\sum_{i=1}^n
		\left[
		\det(v_i,w_{i+1})+
		\det(w_i,v_{i+1})
		\right].
		\]
		In particular,
		\[
		D^2\Acal(c)[v,v]
		=
		\sum_{i=1}^n\det(v_i,v_{i+1}).
		\]
	\end{theorem}
	
	\begin{proof}
		Expand
		\[
		\Acal(c+sv+tw)
		=
		\frac12\sum_{i=1}^n
		\det(c_i+sv_i+tw_i,c_{i+1}+sv_{i+1}+tw_{i+1}).
		\]
		The coefficient of \(st\) in the \(i\)-th determinant is
		\[
		\det(v_i,w_{i+1})+\det(w_i,v_{i+1}),
		\]
		which gives the bilinear formula. Setting \(w=v\) yields the pure
		second variation.
	\end{proof}
	
	\begin{definition}
		Suppose \(c\) satisfies the first order multiplier equations for
        a chosen multiplier vector
		\((\lambda_{ij})_{\{i,j\}\in E(c)}\). The associated constrained
		second variation is
		\[
		Q_c(v)
		=
		D^2\Acal(c)[v,v]
		-
		2\sum_{\{i,j\}\in E(c)}
		\lambda_{ij}\lVert v_i-v_j\rVert ^2,
		\qquad
		v\in(\R^2)^n,
		\]
		and we write \(B_c\) for its polar bilinear form. For nonnegative
		multipliers, the equality space determined by the positive multiplier
		edges may be larger than \(\Ical(c)\). The ambient definition of
		\(Q_c\) allows both restrictions.
	\end{definition}
	
	\begin{theorem}\label{thm:constrained-second-variation}
		Suppose \(c\) satisfies the multiplier identity of
		Proposition~\ref{prop:first-order-no-ascent} for a fixed multiplier
		vector \(\lambda\). Let \(x(t)\) be a \(C^2\) curve satisfying
		\[
		g_{ij}(x(t))=0
		\qquad
		(\{i,j\}\in E(c),\ \lambda_{ij}\neq 0)
		\]
		for all sufficiently small \(t\), with \(x(0)=c\) and \(x'(0)=v\).
		Then
		\[
		\frac{d^2}{dt^2}\bigg|_{t=0}\Acal(x(t))=Q_c(v).
		\]
	\end{theorem}
	
	\begin{proof}
		Set \(a=x''(0)\), so that
		\[
		\frac{d^2}{dt^2}\bigg|_{t=0}\Acal(x(t))
		=
		D^2\Acal(c)[v,v]+D\Acal(c)[a].
		\]
		Differentiating \(g_{ij}(x(t))=0\) twice at \(t=0\) gives, for every
		active edge \(\{i,j\}\) with nonzero multiplier,
		\[
		D^2g_{ij}(c)[v,v]+Dg_{ij}(c)[a]=0,
		\]
		while the multiplier identity gives
		\[
		D\Acal(c)[a]
		=
		\sum_{\{i,j\}\in E(c)}
		\lambda_{ij}Dg_{ij}(c)[a].
		\]
		Combining the two and noting that the remaining multipliers vanish,
		\[
		D\Acal(c)[a]
		=
		-\sum_{\{i,j\}\in E(c)}
		\lambda_{ij}D^2g_{ij}(c)[v,v].
		\]
		Since \(D^2g_{ij}(c)[v,v]=2\lVert v_i-v_j\rVert ^2\), we conclude that
		\[
		\frac{d^2}{dt^2}\bigg|_{t=0}\Acal(x(t))
		=
		D^2\Acal(c)[v,v]
		-
		2\sum_{\{i,j\}\in E(c)}
		\lambda_{ij}\lVert v_i-v_j \rVert ^2
		=
		Q_c(v).
		\]
	\end{proof}
	
	\begin{lemma}
		\label{lem:rigid-kernel}
		Suppose \(c\) satisfies the multiplier identity with multiplier vector
		\(\lambda\). Then
		\[
		B_c(r,w)=0
		\]
		for every \(r\in R(c)\) and every \(w\in(\R^2)^n\).
	\end{lemma}
	
	\begin{proof}
		For translations the identity is immediate, and for rotations it
		suffices, by linearity, to take
		\[
		r_i=Jc_i.
		\]
		The map \(e^{tJ}\) acts componentwise on \((\R^2)^n\), so rotational
		invariance of the area gives
		\[
		\Acal(e^{tJ}(c+sw))=\Acal(c+sw).
		\]
		Differentiating first in \(s\) and then in \(t\) at \(s=t=0\) yields
		\[
		D^2\Acal(c)[r,w]+D\Acal(c)[Jw]=0,
		\]
		while the multiplier identity gives
		\[
		D\Acal(c)[Jw]
		=
		2\sum_{\{i,j\}\in E(c)}
		\lambda_{ij}\langle c_i-c_j,Jw_i-Jw_j\rangle.
		\]
		Since \(J\) is skew adjoint and \(r_i-r_j=J(c_i-c_j)\),
		\[
		\langle c_i-c_j,Jw_i-Jw_j\rangle
		=
		-\langle J(c_i-c_j),w_i-w_j\rangle
		=
		-\langle r_i-r_j,w_i-w_j\rangle,
		\]
		so that
		\[
		D\Acal(c)[Jw]
		=
		-2\sum_{\{i,j\}\in E(c)}
		\lambda_{ij}\langle r_i-r_j,w_i-w_j\rangle
		\]
		and therefore
		\[
		D^2\Acal(c)[r,w]
		=
		2\sum_{\{i,j\}\in E(c)}
		\lambda_{ij}\langle r_i-r_j,w_i-w_j\rangle,
		\]
		which is precisely the assertion \(B_c(r,w)=0\).
	\end{proof}

Relative to the Euclidean inner product on \((\R^2)^n\), the restriction
of \(Q_c\) to \(\Ical(c)\cap R(c)^\perp\) has a spectrum independent of
the chosen orthonormal basis and invariant under orientation preserving
Euclidean motions. We call it the reduced spectrum of \(c\) for the chosen
multipliers.

	\begin{proposition}\label{prop:second-order-exclusion}
		Let \(c\) be a positively oriented small \(n\)-gon. Suppose there exist
		multipliers
		\[
		\lambda_{ij}>0
		\qquad
		(\{i,j\}\in E(c))
		\]
		such that
		\[
		\nabla_{c_i}\Acal(c)
		=
		2\sum_{j:\{i,j\}\in E(c)}
		\lambda_{ij}(c_i-c_j)
		\qquad
		(i=1,\ldots,n),
		\]
		and suppose \(Q_c\) is negative definite on
		\(\Ical(c)\cap R(c)^\perp\). Then every \(C^2\) feasible curve
		\(x(t)\), \(0\leq t<\varepsilon\), of small \(n\)-gons with
		\(x(0)=c\), fixed cyclic labelling and initial velocity
		\[
		x'(0)\notin R(c)
		\]
		satisfies
		\[
		\liminf_{t\downarrow 0}
		\frac{\Acal(x(t))-\Acal(c)}{t^2}<0.
		\]
	\end{proposition}
	
	\begin{proof}
		Let \(x(t)\), \(0\leq t<\varepsilon\), be a \(C^2\) feasible curve with
		\(x(0)=c\), and set
		\[
		v=x'(0),\qquad a=x''(0).
		\]
		Feasibility gives
		\[
		\lVert x_i(t)-x_j(t) \rVert\leq 1
		\]
		for every \(\{i,j\}\in E(c)\). Since \(\lVert c_i-c_j\rVert=1\),
		differentiating at \(t=0\) gives
		\[
		\langle c_j-c_i,v_j-v_i\rangle\leq 0,
		\]
		so \(v\in\Kcal(c)\) and
		Proposition~\ref{prop:first-order-no-ascent} gives
		\(D\Acal(c)[v]\leq 0\). If \(D\Acal(c)[v]<0\), then
		\[
		\Acal(x(t))-\Acal(c)
		=
		tD\Acal(c)[v]+O(t^2),
		\]
		whence
		\[
		\liminf_{t\downarrow 0}
		\frac{\Acal(x(t))-\Acal(c)}{t^2}
		=
		-\infty<0.
		\]
		Assume therefore that \(D\Acal(c)[v]=0\). Because all active
		multipliers are positive, Lemma~\ref{lem:positive-multipliers} forces
		\[
		\langle c_j-c_i,v_j-v_i\rangle=0
		\qquad
		(\{i,j\}\in E(c)),
		\]
		so \(v\in\Ical(c)\). Write
		\[
		v=r+z,
		\qquad
		r\in R(c),
		\qquad
		z\in\Ical(c)\cap R(c)^\perp,
		\]
		and note that \(z\neq 0\), since \(v\notin R(c)\).
		
		For each active edge, the relations \(g_{ij}(x(t))\leq 0\),
		\(g_{ij}(c)=0\) and \(Dg_{ij}(c)[v]=0\) give, by Taylor expansion,
		\[
		D^2g_{ij}(c)[v,v]+Dg_{ij}(c)[a]\leq 0.
		\]
		The multiplier identity gives
		\[
		D\Acal(c)[a]
		=
		\sum_{\{i,j\}\in E(c)}
		\lambda_{ij}Dg_{ij}(c)[a],
		\]
		hence
		\[
		D\Acal(c)[a]
		\leq
		-\sum_{\{i,j\}\in E(c)}
		\lambda_{ij}D^2g_{ij}(c)[v,v],
		\]
		and, since \(D^2g_{ij}(c)[v,v]=2\lVert v_i-v_j\rVert^2\),
		\[
		\frac{d^2}{dt^2}\bigg|_{t=0}\Acal(x(t))
		\leq
		Q_c(v).
		\]
		By Lemma~\ref{lem:rigid-kernel},
		\[
		Q_c(v)=Q_c(r+z)=Q_c(z),
		\]
		and \(Q_c(z)<0\) by negative definiteness, because \(z\neq 0\).
		As \(D\Acal(c)[v]=0\), Taylor expansion gives
		\[
		\Acal(x(t))-\Acal(c)
		=
		\frac12t^2\frac{d^2}{dt^2}\bigg|_{t=0}\Acal(x(t))
		+
		o(t^2),
		\]
		so that
		\[
		\liminf_{t\downarrow 0}
		\frac{\Acal(x(t))-\Acal(c)}{t^2}
		\leq
		\frac12Q_c(v)
		=
		\frac12Q_c(z)<0. \qedhere
		\]
	\end{proof}
	
	\begin{proposition}
		\label{prop:slice-local-max}
		Let \(c\) be a positively oriented small \(n\)-gon, and assume:
		\begin{enumerate}[label=\textup{(\roman*)}]
			\item The active equality set
			\[
			M=\{x:g_{ij}(x)=0\text{ for all }\{i,j\}\in E(c)\}
			\]
			is a \(C^2\) submanifold near \(c\), with tangent space
			\(T_cM=\Ical(c)\).
			
			\item There exist real multipliers \(\lambda_{ij}\)
			satisfying the multiplier identity
			\[
			\nabla_{c_i}\Acal(c)
			=
			2\sum_{j:\{i,j\}\in E(c)}
			\lambda_{ij}(c_i-c_j)
			\qquad
			(i=1,\ldots,n).
			\]
			
			\item \(Q_c\) is negative definite on
			\(\Ical(c)\cap R(c)^\perp\).
		\end{enumerate}
		Then \(c\) is a strict local maximiser of \(\Acal\) on \(M\), modulo
		Euclidean motions: there is a neighbourhood \(U\) of \(c\) in \(M\)
		such that every \(x\in U\) satisfies
		\[
		\Acal(x)\leq \Acal(c),
		\]
		with equality only if \(x\) lies in the Euclidean orbit of \(c\).
		The same conclusion holds in the local fixed active diameter class.
	\end{proposition}
	
	\begin{proof}
		Since \(c\) is positively oriented, its orbit under orientation
		preserving Euclidean motions is a three dimensional embedded
		submanifold with tangent space \(R(c)\).
		The inclusion \(R(c)\subseteq T_cM=\Ical(c)\) implies that the
		affine space \(c+R(c)^\perp\) meets \(M\) transversely at \(c\).
		The implicit function theorem therefore makes
		\(\mathcal S=M\cap(c+R(c)^\perp)\), locally at \(c\), a \(C^2\)
		submanifold with
		\[
		T_c\mathcal S=\Ical(c)\cap R(c)^\perp.
		\]
		Identify the tangent space of the motion group at the identity with
		\(R(c)\) through its action at \(c\). The motion action on the slice
		has differential \((r,z)\mapsto r+z\) at the identity and \(c\).
		This is an isomorphism from \(R(c)\times T_c\mathcal S\) onto
		\(T_cM\). The inverse function theorem gives a local \(C^2\)
		parametrisation of \(M\) near \(c\).
		
		The multiplier identity makes \(D\Acal(c)\) vanish on \(T_cM\).
		By Theorem~\ref{thm:constrained-second-variation}, the Hessian of
		\(\Acal|_{\mathcal S}\) at \(c\) is the restriction of \(Q_c\) to
		\(T_c\mathcal S\), which is negative definite by (iii). The second
		derivative test gives a strict local maximum on the slice when it
		has positive dimension. A zero-dimensional slice consists locally
		only of \(c\), with the same conclusion.
		
		Area is invariant under the motion action, so the inequality and its
		equality condition extend to \(M\). Since
		\(\delta(c)>0\) by definition, all inactive distances
		remain less than \(1\) in a sufficiently small neighbourhood.
		There \(M\) coincides with the fixed active diameter class.
	\end{proof}

\begin{definition}
For \(x,c\in(\R^2)^n\), define
\[
\operatorname{dist}(x,c)=
\min_{\theta\in[0,2\pi],\ a\in\R^2}
\Bigl\lVert
\bigl(e^{\theta J}x_i+a-c_i\bigr)_{i=1}^n
\Bigr\rVert.
\]

\end{definition}

The minimum is attained. For fixed \(\theta\), the
optimal translation is \(a=\bar c-e^{\theta J}\bar x\),
where bars denote centroids; the remaining function
of \(\theta\) is continuous on \([0,2\pi]\).

For an \(m\times N\) matrix \(B\), with \(1\leq m\leq N\),
let \(s_{\min}(B)\) denote its smallest singular
value, the square root of the smallest eigenvalue
of \(BB^\top\).

	\begin{theorem}
	    \label{thm:quant-stability}
		Let \(c\) be a small \(n\)-gon with \(1\leq |E(c)|\leq 2n\), and let
		\(\Delta\) be the maximum degree of \(G(c)\). Assume:
		\begin{enumerate}[label=\textup{(\roman*)}]
			\item \(s_{\min}(Dg(c))\geq s_L>0\).
			
			\item There are multipliers satisfying the multiplier identity of
			Proposition~\ref{prop:first-order-no-ascent} with
			\(0<\lambda_L\leq\lambda_{ij}\leq\lambda_U\) for every
			\(\{i,j\}\in E(c)\).
			
			\item \(Q_c(v)\leq-\mu\lVert v\rVert^2\) for every
			\(v\in\Ical(c)\cap R(c)^\perp\), where \(\mu>0\).
		\end{enumerate}
		Set
		\[
		q=1+4\Delta\lambda_U,
		\qquad K=\frac q2+\frac{q^2}{\mu}+\frac{\mu}{4},
		\]
		and
		\[
		\rho=\frac{s_L}{K}\min\Bigl\{\lambda_L,\frac{\mu}{16\Delta}\Bigr\},
		\qquad \kappa=\frac{\mu}{8}.
		\]
		Then every \(x\in(\R^2)^n\) satisfying
		\[
		g_{ij}(x)\leq0\quad(\{i,j\}\in E(c)),
		\qquad \operatorname{dist}(x,c)\leq\rho
		\]
		obeys
		\[
		\Acal(x)\leq\Acal(c)-\kappa\operatorname{dist}(x,c)^2.
		\]
		In particular, this holds for every small \(n\)-gon within distance
		\(\rho\) of \(c\). Equality \(\Acal(x)=\Acal(c)\) forces such an \(x\)
		to lie in the Euclidean orbit of \(c\).
\end{theorem}
	
	\begin{proof}
		Choose an orientation preserving Euclidean motion attaining
		\(\operatorname{dist}(x,c)\), and let \(y\) be its image of \(x\).
		Set \(h=y-c\) and
		\(r=\lVert h\rVert=\operatorname{dist}(x,c)\).
		Minimality under translations and rotations gives, for
		\(b\in\R^2\) and \(\omega\in\R\),
		\[
		0=\frac{d}{ds}\bigg|_{s=0}
		\sum_{i=1}^n\bigl\lVert e^{s\omega J}y_i+sb-c_i\bigr\rVert^2
		=2\sum_{i=1}^n\langle\omega Jy_i+b,h_i\rangle.
		\]
		Taking \(\omega=0\) gives \(\sum_i h_i=0\). Taking \(b=0\) and
		\(\omega=1\), then using \(y_i=c_i+h_i\) and
		\(\langle Jh_i,h_i\rangle=0\), gives
		\(\sum_i\langle Jc_i,h_i\rangle=0\).
		Thus \(h\in R(c)^\perp\). The motion preserves area and pairwise
		distances, so it suffices to consider \(c+h\) with \(r\leq\rho\)
		and \(g_{ij}(c+h)\leq0\) on every edge of \(E(c)\).
		
		Consider the quadratic Lagrangian
		\[
		L(u)=\Acal(u)-\sum_{\{i,j\}\in E(c)}\lambda_{ij}g_{ij}(u).
		\]
		The multiplier identity gives \(DL(c)=0\); also
		\(L(c)=\Acal(c)\) and \(D^2L(c)[h,h]=Q_c(h)\).
		Since \(L\) is quadratic, its expansion is exact:
		\[
		\Acal(c+h)-\Acal(c)
		=\sum_{\{i,j\}\in E(c)}\lambda_{ij}g_{ij}(c+h)
		+\frac12Q_c(h).
		\]
		For each active edge, the constraint expansion is likewise exact:
		\[
		g_{ij}(c+h)=Dg_{ij}(c)[h]+\lVert h_i-h_j\rVert^2.
		\]
		Write the active constraint slacks and their sum as
		\[
		\sigma_{ij}=-g_{ij}(c+h)\geq0,
		\qquad S=\sum_{\{i,j\}\in E(c)}\sigma_{ij}.
		\]
		The multiplier bound in (ii) gives
		\[
		\Acal(c+h)-\Acal(c)\leq-\lambda_LS+\frac12Q_c(h).
		\]
		
		The graph degree bound, valid for every \(u\in(\R^2)^n\), is
		\[
		\sum_{\{i,j\}\in E(c)}\lVert u_i-u_j\rVert^2
		\leq2\sum_{i=1}^n\deg_{G(c)}(i)\lVert u_i\rVert^2
		\leq2\Delta\lVert u\rVert^2.
		\]
		Decompose \(h=t+z\) orthogonally, with \(t\in\ker Dg(c)\) and
		\(z\in(\ker Dg(c))^\perp\). Since \(R(c)\subseteq\ker Dg(c)\) and
		\(h\in R(c)^\perp\), we have
		\(t\in\Ical(c)\cap R(c)^\perp\). Hypothesis (i) gives
		\[
		\lVert Dg(c)z\rVert\geq s_L\lVert z\rVert,
		\qquad \lVert t\rVert^2+\lVert z\rVert^2=r^2.
		\]
		The constraint expansion gives, on every active edge,
		\[
		Dg_{ij}(c)[z]=-\sigma_{ij}-\lVert h_i-h_j\rVert^2.
		\]
		Taking norms and applying the degree bound to \(h\) yields
		\[
		\begin{gathered}
			\lVert Dg(c)z\rVert
			\leq\sum_{\{i,j\}\in E(c)}
			\bigl(\sigma_{ij}+\lVert h_i-h_j\rVert^2\bigr)
			\leq S+2\Delta r^2,\\
			\lVert z\rVert\leq\frac{S+2\Delta r^2}{s_L}.
		\end{gathered}
		\]
		
		For every \(v\in(\R^2)^n\), Theorem~\ref{thm:second-variation} and the degree
		bound give
		\[
		\begin{gathered}
			|D^2\Acal(c)[v,v]|
			\leq\sum_{i=1}^n\lVert v_i\rVert\,\lVert v_{i+1}\rVert
			\leq\lVert v\rVert^2,\\
			0\leq 2\sum_{\{i,j\}\in E(c)}\lambda_{ij}\lVert v_i-v_j\rVert^2
			\leq4\Delta\lambda_U\lVert v\rVert^2.
		\end{gathered}
		\]
		Thus \(|Q_c(v)|\leq q\lVert v\rVert^2\). The symmetric operator
		associated with \(B_c\) has norm at most \(q\), so
		\(|B_c(u,v)|\leq q\lVert u\rVert\lVert v\rVert\).
		By (iii) and Young's inequality,
		\[
		\begin{gathered}
			Q_c(h)=Q_c(t)+2B_c(t,z)+Q_c(z),\\
			Q_c(h)\leq-\mu\lVert t\rVert^2
			+2q\lVert t\rVert\lVert z\rVert+q\lVert z\rVert^2,\\
			Q_c(h)\leq-\frac{\mu}{2}\lVert t\rVert^2
			+\Bigl(q+\frac{2q^2}{\mu}\Bigr)\lVert z\rVert^2.
		\end{gathered}
		\]
		Since \(\lVert t\rVert^2=r^2-\lVert z\rVert^2\), the definition
		of \(K\) gives
		\[
		\frac12Q_c(h)\leq-\frac{\mu}{4}r^2+K\lVert z\rVert^2.
		\]
		
		The estimate for \(\lVert z\rVert\), together with
		\(\lVert z\rVert\leq r\), gives
		\[
		\lVert z\rVert^2\leq\frac{r}{s_L}(S+2\Delta r^2).
		\]
		Substitution into the area estimate yields
		\[
		\Acal(c+h)-\Acal(c)
		\leq-\Bigl(\lambda_L-\frac{Kr}{s_L}\Bigr)S
		-r^2\Bigl(\frac{\mu}{4}-\frac{2\Delta Kr}{s_L}\Bigr).
		\]
		If \(r\leq\rho\), then \(Kr/s_L\leq\lambda_L\) and
		\(2\Delta Kr/s_L\leq\mu/8\). Since \(S\geq0\), we obtain
		\[
		\Acal(c+h)-\Acal(c)\leq-\frac{\mu}{8}r^2=-\kappa r^2.
		\]
		Equality \(\Acal(x)=\Acal(c)\) forces \(r=0\), so \(x\) lies in the
		Euclidean orbit of \(c\).
	\end{proof}
	
	\begin{samepage}
		\begin{theorem}
			\label{thm:ambient-local-max}
			Let \(c\) be a positively oriented small \(n\)-gon with
			\(E(c)\neq\emptyset\). Assume:
			\begin{enumerate}[label=\textup{(\roman*)}]
				\item The active constraint matrix \(Dg(c)\) has full row rank.
				
				\item There exist multipliers \(\lambda_{ij}>0\) satisfying
				\[
				\nabla_{c_i}\Acal(c)
				=
				2\sum_{j:\{i,j\}\in E(c)}
				\lambda_{ij}(c_i-c_j)
				\qquad
				(i=1,\ldots,n).
				\]
				
				\item \(Q_c\) is negative definite on
				\[
				\Ical(c)\cap R(c)^\perp.
				\]
			\end{enumerate}
			Then \(c\) is a strict local maximiser of \(\Acal\), modulo
			orientation preserving Euclidean motions, among all small
			\(n\)-gons. Thus there is a neighbourhood \(U\) of \(c\) such
			that every small \(n\)-gon \(x\in U\) satisfies
			\[
			\Acal(x)\leq\Acal(c),
			\]
			with equality only if \(x\) lies in the Euclidean orbit of \(c\).
			The neighbourhood can also be chosen so that
			\(E(x)\subseteq E(c)\) for every such \(x\).
		\end{theorem}
	\end{samepage}
	
	\begin{proof}
		By (i), \(s_{\min}(Dg(c))>0\). The finitely many multipliers in (ii)
		lie in an interval \([\lambda_L,\lambda_U]\) with \(\lambda_L>0\).
		Set \(V=\Ical(c)\cap R(c)^\perp\).
		If \(V=\{0\}\), choose any \(\mu>0\). Otherwise, (iii) and compactness
		of the unit sphere of \(V\) give \(\mu>0\). In either case,
		\(Q_c(v)\leq-\mu\lVert v\rVert^2\) for all \(v\in V\).
		
		Every small \(n\)-gon \(x\) satisfies \(\lVert x_i-x_j\rVert\leq1\) on
		\(E(c)\). Theorem~\ref{thm:quant-stability} therefore gives
		\(\rho>0\) such that
		\[
		\Acal(x)\leq\Acal(c)-\frac{\mu}{8}\operatorname{dist}(x,c)^2
		\qquad\text{if }\operatorname{dist}(x,c)\leq\rho.
		\]
		Let \(U=\{x:\lVert x-c\rVert<\rho\}\). For every small \(n\)-gon
		\(x\in U\), we have
		\(\operatorname{dist}(x,c)\leq\lVert x-c\rVert<\rho\).
		Hence \(\Acal(x)\leq\Acal(c)\), with equality only when \(x\) lies
		in the Euclidean orbit of \(c\).
		
		The inactive margin \(\delta(c)\) is positive by definition.
		Shrinking \(U\) so that
		\(\max_i\lVert x_i-c_i\rVert<\delta(c)/2\)
		keeps every inactive distance strictly below \(1\), and therefore
		\(E(x)\subseteq E(c)\) for every small \(n\)-gon \(x\in U\).
	\end{proof}
	
	\begin{remark}
		\label{rem:exactness}
		For a Graham type diameter graph, \(\Delta=3\) in
		Theorem~\ref{thm:quant-stability}: the cycle vertex incident with
		the pendant edge has degree \(3\), the pendant vertex has degree \(1\),
		and every other vertex has degree \(2\).
	\end{remark}

	\section{The quadrilateral case}
	\label{sec:models}
	The quadrilateral \(Q_4\) illustrates how local maximality can coexist
	with vanishing active multipliers and a degenerate second
	variation. Its explicit description also shows how the remaining
	diameter inequalities restrict the directions in the reduced kernel.
	
	\begin{proposition}
		\label{prop:q4}
		Let \(a=1-\sqrt3/2\) and define \(Q_4=c=(c_1,c_2,c_3,c_4)\) by
		\[
		c_1=(0,0),\qquad
		c_2=\left(\frac12,a\right),\qquad
		c_3=(0,1),\qquad
		c_4=\left(-\frac12,a\right).
		\]
		Then \(\Acal(c)=1/2\). The active diameter edges are
		\[
		E(c)=\bigl\{\{1,3\},\{2,3\},\{2,4\},\{3,4\}\bigr\}.
		\]
		The first order multiplier equations are satisfied by
		\[
		\lambda_{13}=\frac14,\qquad
		\lambda_{24}=\frac14,\qquad
		\lambda_{23}=0,\qquad
		\lambda_{34}=0.
		\]
		Thus \(Q_4\) satisfies the first order no ascent criterion. Let
		\[
		\Ical^+(c)=
		\left\{
		v\in(\R^2)^4:
		\langle c_3-c_1,v_3-v_1\rangle=0,
		\ \langle c_4-c_2,v_4-v_2\rangle=0
		\right\}.
		\]
		Only the positive multiplier edges enter this equality space. The constrained
		second variation on it is
		\[
		Q_c(v)
		=
		D^2\Acal(c)[v,v]
		-\frac12\lVert v_1-v_3\rVert^2
		-\frac12\lVert v_2-v_4\rVert^2.
		\]
		Moreover,
		\[
		\dim \Ical^+(c)=6,\qquad
		\dim R(c)=3,\qquad
		\dim(\Ical^+(c)\cap R(c)^\perp)=3.
		\]
		On \(\Ical^+(c)\cap R(c)^\perp\) the reduced eigenvalues are
		\[
			\{-2,\ 0,\ 0\}.
		\]
		Within the zero eigenspace, the full diameter inequalities restrict feasible
		directions to a proper cone.
	\end{proposition}
	
	\begin{proof}
		The shoelace formula gives \(\Acal(c)=1/2\). For any small quadrilateral, the absolute signed area is at most half the product of its diagonal lengths, each at most \(1\). Thus its area cannot exceed \(1/2\). The active diameter equalities
		follow from
		\[
		\lVert c_1-c_3\rVert=1,\qquad
		\lVert c_2-c_3\rVert=1,\qquad
		\lVert c_2-c_4\rVert=1,\qquad
		\lVert c_3-c_4\rVert=1.
		\]
		The remaining two distances are
		\[
		\lVert c_1-c_2\rVert=\lVert c_1-c_4\rVert
		=
		\sqrt{\frac14+\left(1-\frac{\sqrt3}{2}\right)^2}<1.
		\]
		By Theorem~\ref{thm:first-variation}, the gradients are
		\[
		\nabla_{c_1}\Acal=\left(0,-\frac12\right),\quad
		\nabla_{c_2}\Acal=\left(\frac12,0\right),\quad
		\nabla_{c_3}\Acal=\left(0,\frac12\right),\quad
		\nabla_{c_4}\Acal=\left(-\frac12,0\right).
		\]
		Substitution gives the multiplier equations with the displayed multipliers.
		The two equations defining \(\Ical^+(c)\) are
		\[
		(v_3-v_1)_2=0,\qquad
		(v_4-v_2)_1=0.
		\]
		Hence \(\dim\Ical^+(c)=8-2=6\). Since \(R(c)\) has dimension \(3\) and is
		contained in \(\Ical^+(c)\), the reduced space
		\(\Ical^+(c)\cap R(c)^\perp\) has dimension \(3\). A basis of
		\(\Ical^+(c)\cap R(c)^\perp\) is given by
		\[
		u_1=(1,1,0,-2,-1,1,0,0),
		\]
		\[
		u_2=(-\sqrt3,0,1,0,-2+\sqrt3,0,1,0),
		\]
		\[
		u_3=(1,0,0,-1,-1,0,0,1),
		\]
		where the coordinates are ordered as
		\[
		(v_{1x},v_{1y},v_{2x},v_{2y},v_{3x},v_{3y},v_{4x},v_{4y}).
		\]
		With respect to this basis, the matrix of \(Q_c\) and the Euclidean Gram
		matrix are
		\[
		Q=
		\begin{pmatrix}
			-8 & -4+4\sqrt3 & -8\\
			-4+4\sqrt3 & -8+4\sqrt3 & -4+4\sqrt3\\
			-8 & -4+4\sqrt3 & -8
		\end{pmatrix},
		\]
		and
		\[
		G=
		\begin{pmatrix}
			8 & 2-2\sqrt3 & 4\\
			2-2\sqrt3 & 12-4\sqrt3 & 2-2\sqrt3\\
			4 & 2-2\sqrt3 & 4
		\end{pmatrix}.
		\]
		The reduced eigenvalues are the roots of the generalised characteristic
		equation
		\[
		\det(Q-\nu G)=0.
		\]
		A direct calculation gives
		\[
		\det(Q-\nu G)
		=
		32\nu^2(\sqrt3-4)(\nu+2).
		\]
		Since \(\sqrt3-4\neq 0\), the reduced spectrum is exactly
		\[
			\{-2,\ 0,\ 0\}.
		\]
		To identify the kernel, keep \(c_1=(0,0)\), \(c_3=(0,1)\), and put
		\(c_2=(t+1/2,b)\), \(c_4=(t-1/2,b)\). The two diagonals remain orthogonal
		and have unit length, so the area stays \(1/2\). The area and the two
		positive-multiplier constraints are constant along this family, so
		\(Q_c\) vanishes on its tangent velocities. These velocities are
		independent modulo rigid motions and, after orthogonal projection onto \(R(c)^\perp\), span the two-dimensional reduced
		kernel. Near \((t,b)=(0,a)\), the remaining active inequalities are
\[
(t+1/2)^2+(b-1)^2\leq1,
\qquad
(t-1/2)^2+(b-1)^2\leq1.
\]
Their first order inequalities at \((0,a)\) are
\[
\dot t-\sqrt3\,\dot b\leq0,
\qquad
-\dot t-\sqrt3\,\dot b\leq0.
\]
Thus feasible velocities satisfy
\(\dot b\geq|\dot t|/\sqrt3\), and \(Q_4\) is a corner of the
feasible part of this optimal family.
	\end{proof}

	\section{Spectral Verification }
\label{sec:lspectra}
	
	LSPectra~\cite{lspectra} encloses the reduced spectrum and verifies the
	hypotheses of Theorem~\ref{thm:quant-stability} and
	Theorem~\ref{thm:ambient-local-max} for each even \(n\) with
	\(6\le n\le\nmax\). The program uses outward rounded interval arithmetic
	to enclose a stationary configuration and bound the quantities entering
	the stability estimate. These bounds hold at the exact stationary configuration and are given in Table~\ref{tab:records}. Employing these bounds and \(\Delta=3\), the constants \(\rho_n\) and \(\kappa_n\) of
	Theorem~\ref{thm:quant-stability} are evaluated in exact rational arithmetic and then rounded downward. The resulting values of \(\rho_n\) and \(\kappa_n\) appear in Section~\ref{sec:main}.
	
	For each even \(n\) with \(6\le n\le\nmax\), put \(h=n/2\).
    The \(m=n\) pairs in \(\mathcal E\) are listed in the order given below.
    The input declares the pendant pair \(\{1,h+1\}\) and the cycle edges
    \(\{k,k+h-1\}\), \(\{k,k+h\}\) for \(2\le k\le h\), and \(\{h+1,n\}\).
    We write \(\mathcal E\) for this list and also for its underlying set of pairs.

\begin{remark}
For a strictly convex realisation of the cycle-plus-pendant
diameter graph, the active constraint gradients are linearly
independent. Indeed, in any linear relation among these gradients,
the block corresponding to the pendant vertex forces the coefficient
of its sole incident edge to vanish. At each cycle vertex, the
remaining two terms are scalar multiples of the two incident cycle
chord vectors. These vectors are nonparallel by strict convexity,
so both coefficients vanish.
\end{remark}

    Define
	\[
	g(z)=\bigl(\lVert z_i-z_j\rVert^2-1\bigr)_{\{i,j\}\in\mathcal E}
	\qquad (z\in(\R^2)^n).
	\]
	The same pairs define \(g\) and its Jacobian \(Dg\) throughout the
	verification box. The input coordinates for \(n=6,8,10\) come from
	Mossinghoff~\cite[Table~4]{mossinghoff}; for \(n\ge12\), they are generated
	from OPTIGON angles~\cite{optigon}. The files accompany LSPectra and
	provide starting points for the verification.
	
	To solve the stationarity equations, we fix a gauge to remove
    the freedom associated with Euclidean motions.
    Set \(\Sigma(c)=[\,\mathbf 1_x\ \mathbf 1_y\ Jc\,]\),
    whose columns span \(R(c)\), and impose the gauge
    \(\Gamma c=b\) given by
    \[
    c_1=(0,0),\qquad (c_k)_1=0,\qquad k=\lfloor n/2\rfloor+1.
    \]
    Here \(\mathbf 1_x=((1,0),\ldots,(1,0))\) and
    \(\mathbf 1_y=((0,1),\ldots,(0,1))\) are the translation
    generators, \(\Gamma c=((c_1)_1,(c_1)_2,(c_k)_1)^\top\),
    and \(b=0\in\R^3\).
    The bordered system~\cite{govaerts2000} is
	\[
	\widehat F(c,\lambda,\xi)=
	\begin{pmatrix}
		\nabla\Acal(c)-Dg(c)^\top\lambda+\Sigma(c)\xi\\
		g(c)\\
		\Gamma c-b
	\end{pmatrix},
	\]
	with \((c,\lambda,\xi)\in\R^{2n}\times\R^m\times\R^3\) and
	\(\widehat F(c,\lambda,\xi)\in\R^{2n+m+3}\).
	Every stationary pair satisfying the gauge gives a zero with \(\xi=0\).
	The following lemma gives the converse when the area is positive.
	
	\begin{lemma}
		\label{lem:bordering-consistency}
		Let \((c,\lambda,\xi)\) be an exact zero of \(\widehat F\) with
		\(\Acal(c)>0\). Then \(\xi=0\), and \((c,\lambda)\) satisfies the multiplier
		identity
		\[
		\nabla\Acal(c)=Dg(c)^\top\lambda
		\]
		for the pairs of \(\mathcal E\), without a sign restriction on \(\lambda\).
		Once \(\mathcal E=E(c)\) is shown, this is
		the multiplier identity of Proposition~\ref{prop:first-order-no-ascent}.
	\end{lemma}
	
	\begin{proof}
		The area and the constraints \(g\) are invariant under orientation
		preserving Euclidean motions, so
		\[
		\Sigma(c)^\top\nabla\Acal(c)=0,
		\qquad
		\Sigma(c)^\top Dg(c)^\top=0.
		\]
		Applying \(\Sigma(c)^\top\) to the first block of \(\widehat F=0\) gives
		\[
		\Sigma(c)^\top\Sigma(c)\,\xi=0.
		\]
		Since \(\Acal(c)>0\), the points \(c_i\) are not all equal, so the columns of
		\(\Sigma(c)\) are linearly independent and the Gram matrix
		\(\Sigma(c)^\top\Sigma(c)\) is nonsingular. Hence \(\xi=0\), and the first block
		is the stationarity identity.
	\end{proof}
	
	\subsection{The verification}
	\label{subsec:pipeline}
	
	We first refine an approximate solution of the bordered system using Newton's
    method at 256-bit precision. We round the result to binary64 and use it as 
    the centre of an interval box. A Krawczyk inclusion establishes the existence
    of a unique zero of the bordered system in the box, as required by H1. The 
    remaining checks in Table~\ref{tab:verification-checks} are performed on a 
    contracted enclosure \(Z_n\) containing this zero.
    The interval checks over \(Z_n\) provide the bounds given in
    Section~\ref{sec:main}. These include lower bounds for the inactive margin and
    \(s_{\min}(Dg)\), together with lower and upper bounds for
    the active multipliers. They also include the spectral bound
    \(\mu_n\) and an interval enclosing the area.
    Since \(Z_n\) contains the exact zero, these bounds hold at the
    exact stationary configuration. The threshold for H3--H6 and H8
    is \(\tau=2^{-13}\). Failure to meet this threshold does not by
    itself disprove the underlying strict inequality.
	
	For the spectral check, define
	\[
	H(\lambda)=D^2\Acal-
	\sum_{\{i,j\}\in\mathcal E}\lambda_{ij}D^2g_{ij},
	\qquad
	A(z)=\begin{pmatrix}Dg(z)\\ \Sigma(z)^\top\end{pmatrix}.
	\]
	The Hessians are represented in standard Euclidean coordinates. The test
	concerns the orthogonal compression of \(H(\lambda)\) to
	\(\ker A(z)=\ker Dg(z)\cap R(z)^\perp\). At an exact stationary pair
	with \(E(z)=\mathcal E\), this compression represents the
	restriction of \(Q_z\) to \(\Ical(z)\cap R(z)^\perp\).
	
	\begin{table}[htbp]
		\caption{H1 is checked on the initial box \(Z\).
        The remaining conditions are checked on the contracted box \(Z_n\).}
		\label{tab:verification-checks}
		\centering
		\small
		\renewcommand{\arraystretch}{1.15}
		\begin{tabularx}{0.92\linewidth}{@{}l>{\raggedright\arraybackslash}X@{}}
			\toprule
			Check & Condition\\
			\midrule
			H1 & exactly one zero of the bordered system in the box\\
			H2 & full row rank of \(A(z)\), hence of \(Dg(z)\)\\
			H3 & \(|\det(\Gamma\Sigma(z))|\ge\tau\)\\
			H4 & \(\lambda_{ij}\ge\tau\) for every \(\{i,j\}\in\mathcal E\)\\
			H5 & \(1-\lVert z_i-z_j\rVert\ge\tau\) for every \(\{i,j\}\notin\mathcal E\)\\
			H6 & the largest eigenvalue of the orthogonal compression of
			\(H(\lambda)\) to \(\ker Dg(z)\cap R(z)^\perp\) is below \(-\tau\)\\
			H7 & \(\Acal(z)>0\)\\
			H8 & every non-incident vertex has signed distance \(\ge\tau\)
			to each directed polygon side\\
			\bottomrule
		\end{tabularx}
	\end{table}
	
	In H8, the signed distance from \(z_j\) to the directed side
	\(z_i z_{i+1}\) is
	\[
	\frac{\det(z_{i+1}-z_i,z_j-z_i)}{\lVert z_{i+1}-z_i\rVert}.
	\]
	LSPectra also records a lower bound \(D_n\) for the determinants
	\(\det(z_{i+1}-z_i,z_j-z_i)\) over \(Z_n\), used in
	Remark~\ref{rem:simple-competitors}.\footnote{The bounds \(D_n\) are found in the JSON records under \textit{evidence.convexity.det\_lower\_bound}, in \textit{data/records} for \(n=6,8,\ldots,14\) and \textit{data/records\_series/} for \(n=16,18,\ldots,64\).}
    The signed distance is positive to the left of each directed
    boundary side, so H8 alone implies positive orientation and
    strict convexity. H7 is listed
    separately because Lemma~\ref{lem:bordering-consistency} uses it. At the exact zero supplied
	by H1, \(g(c)=0\) and H5 give \(\mathcal E=E(c)\).
	Thus \(c\) is a small \(n\)-gon and \(\ker Dg(c)=\Ical(c)\), as in
	Section~2. Together with the multiplier identity, H2, H4 and H6 then give
	the rank, positivity and second variation hypotheses of
	Theorem~\ref{thm:ambient-local-max} and, with the bounds reported in
	Section~\ref{sec:main}, those of Theorem~\ref{thm:quant-stability}.
	
	To check H6 and enclose the reduced spectrum, the program counts the
	negative eigenvalues of the shifted bordered matrix
	\[
	\mathcal B(s,z,\lambda)=
	\begin{pmatrix}
		H(\lambda)-sI & A(z)^\top\\
		A(z) & 0
	\end{pmatrix}.
	\]
	When \(A(z)\) has full row rank, subtracting its number of rows from this
	count gives the number of reduced eigenvalues below \(s\)~\cite{gould1985}.
	Approximate eigenvalues guide the choice of shifts, including \(-\tau\)
	and \(0\). For each shift, interval evaluation encloses \(\mathcal B\)
	over the box. With \(P\) an approximate inverse, the bound
	\[
	\|I-P\mathcal B(s,z,\lambda)\|_\infty<1
	\]
	throughout the interval matrix enclosure establishes nonsingularity.
	The inertia is therefore constant on its symmetric members and can be
    determined from one representative. 

The program constructs a symmetric interval enclosure
and takes its entrywise lower endpoints to obtain
a symmetric representative \(M\). It then uses a
matrix \(V\) of approximate eigenvectors of \(M\)
to enclose \(V^\top M V\).

    If all Gershgorin discs exclude zero, the number of discs lying entirely in the
    open left half-plane gives the number of negative eigenvalues of
    \(V^\top M V\). Its nonsingularity implies that \(V\) is invertible,
    so Sylvester's law gives the same inertia for \(M\). If this test is
    inconclusive, the program computes the inertia of \(M\) by symmetric
    elimination in exact rational arithmetic.
    Only shifts for which these checks succeed enter the spectral bounds.
    Differences between the counts at successive shifts give the enclosed
    eigenvalue multiplicities. The enclosures cover the entire reduced
    spectrum when the counts at the extreme shifts are \(0\) and
    \(\dim\ker A(z)\). 
For each \(Q_n\) covered by Theorem~\ref{thm:main}, the shift counts
increase by one across each of \(n-3\) pairwise
disjoint intervals contained in \((-\infty,-\tau)\).
Each interval therefore contains exactly one reduced
eigenvalue, proving that the entire reduced spectrum
is simple and strictly negative.

The following results are formalised in
Lean~4~\cite{lean4} as part of LSPectra~\cite{lspectra}:
Lemma~\ref{lem:bordering-consistency}, the Neumann
nonsingularity bound, constant inertia on interval
families, the bordered inertia identity, the
localisation of eigenvalues between shifts, and the
Gershgorin criterion. 
The interval arithmetic and
Theorem~\ref{thm:quant-stability} are not formalised
in Lean. The logic combining the verification checks
into the final verdict is tested but not formally verified.

	\subsection{Implementation details}
	\label{subsec:implementation}
	
	Let \(u_0=(c_0,\lambda_0,\xi_0)\) be the refined centre,
	\(Z\subset\R^{2n+m+3}\) the box about it, and \(Y\) an approximate
	inverse of \(D\widehat F(u_0)\). The program evaluates
	\[
	\mathcal T(Z)=u_0-Y\,\widehat F(u_0)
	+\bigl(I-Y\,D\widehat F(Z)\bigr)(Z-u_0)
	\]
	in interval arithmetic. Strict inclusion
	\(\mathcal T(Z)\subset\operatorname{int}Z\) implies the existence of a
	unique zero in \(Z\) and nonsingularity of the enclosed
	Jacobians~\cite{moore1977}, \cite[Ch.~5]{neumaier1990}.
	Subsequent checks use the contracted enclosure
	\(Z\cap\mathcal T(Z)\). Every zero of \(\widehat F\) in \(Z\) lies in
	\(\mathcal T(Z)\), so this enclosure, recorded for each \(n\) as \(Z_n\),
	contains the zero, which is unique there.

Each endpoint of \(Z_n\) is written to the JSON record as a decimal number that recovers the original binary64 value when read as binary64. The endpoints defining \(Z_n\) are the exact dyadic values represented by these recovered binary64 numbers.
	
	Pivoted QR at the refined centre selects a square column block
	\((A(z))_B\). On the contracted enclosure, an approximate inverse \(Y_B\)
	is used to verify
	\[
	\|I-Y_B(A(z))_B\|_\infty<1.
	\]
	This proves full row rank\footnote{The upper bound on
    the residual norm is stored in each record under
    \textit{evidence.rank.kappa}.}. Solving for the basic variables with interval
	bounds then encloses a basis of \(\ker Dg(z)\cap R(z)^\perp\).
	For H3, the interval evaluation establishes \(|\det(\Gamma\Sigma(z))| \ge\tau\)
	throughout the box, ensuring transversality of the gauge to Euclidean motion orbits.
    This check is required for LSPectra's \textit{STABLE} verdict, but not separately for the proof of Theorem~\ref{thm:main}.
    
	For the stability estimate, the program applies interval Gaussian
    elimination without pivoting to an enclosure of
    \(Dg(z)Dg(z)^\top-s_{n,L}^2I\) over \(Z_n\). It accepts the candidate
    \(s_{n,L}\) if every diagonal pivot has a strictly positive lower
    endpoint. This proves positive definiteness throughout \(Z_n\) and
    hence \(s_{\min}(Dg(z))>s_{n,L}\).
	
	The implementation uses Julia~\cite{julia2017} and
	IntervalArithmetic.jl~\cite{intervalarithmetic}. The verification scripts
	check outward rounding before computation. The software includes input
	files, records and reproduction scripts, together with regression and
	mutation tests. The regression tests recompute the cases and compare
	them with the stored records; the extended series is enabled separately.
	
\section{The main theorem and the verified spectra}
\label{sec:main}

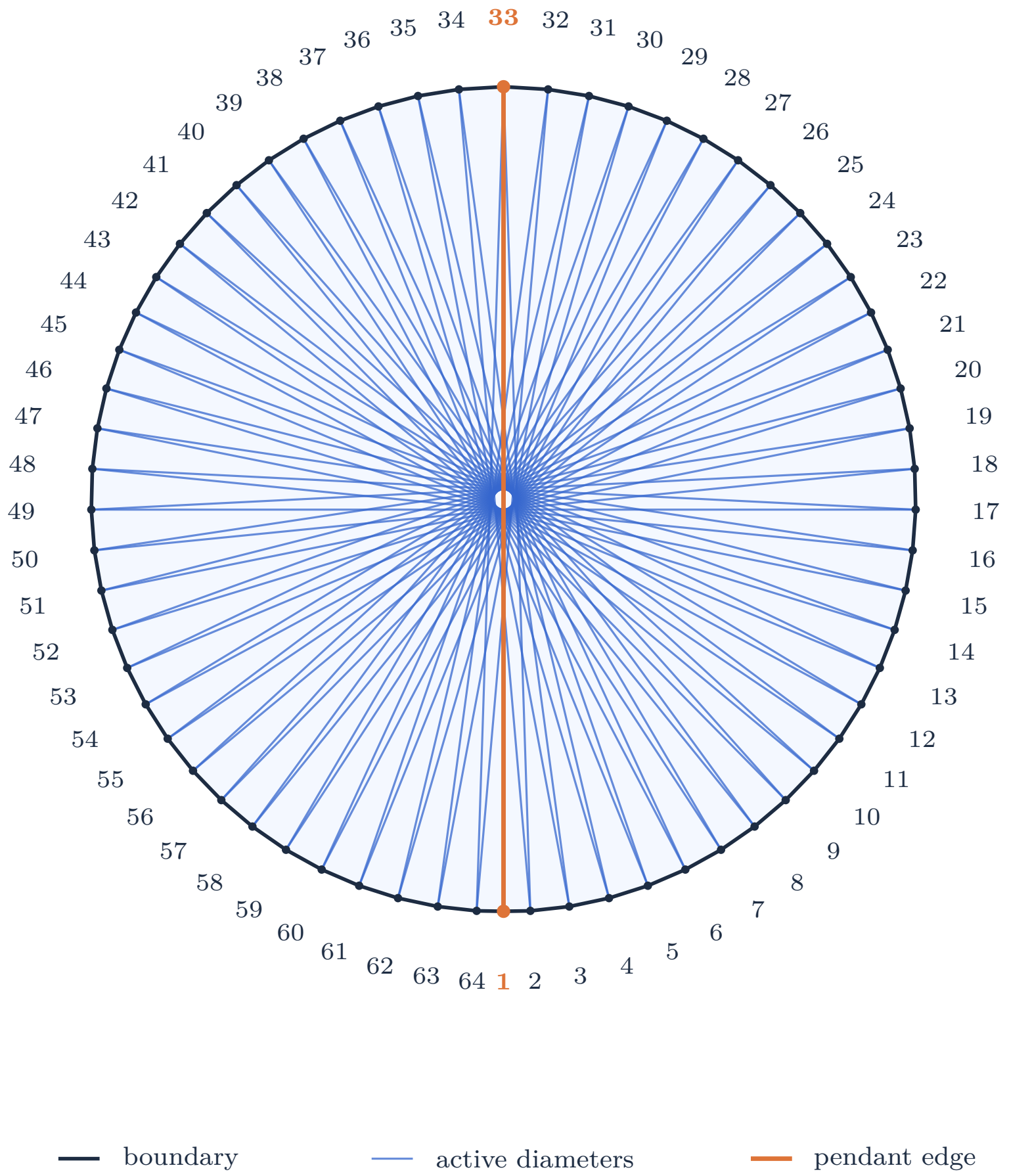
\begin{figure}[!ht]
    \centering
    \definecolor{q64navy}{HTML}{1D2C42}
    \definecolor{q64blue}{HTML}{3667CE}
    \definecolor{q64orange}{HTML}{DE7438}
    \definecolor{q64pale}{HTML}{F4F8FF}
    \begin{tikzpicture}[
        x=8.3cm,y=8.3cm,
        diam/.style={draw=q64blue,line width=0.60pt,opacity=0.75},
        borderedge/.style={draw=q64navy,line width=1.05pt},
        pendant/.style={draw=q64orange,line width=1.30pt},
        dot/.style={circle,fill=q64navy,inner sep=0.85pt},
        endpoint/.style={circle,fill=q64orange,inner sep=1.35pt}
    ]


        \coordinate (G1) at (0.0000000000,0.0000000000);
        \coordinate (G2) at (0.0325860417,0.0005310661);
        \coordinate (G3) at (0.0798158765,0.0057745326);
        \coordinate (G4) at (0.1279693828,0.0159973810);
        \coordinate (G5) at (0.1751078600,0.0310019184);
        \coordinate (G6) at (0.2205644238,0.0506070612);
        \coordinate (G7) at (0.2638644679,0.0746128255);
        \coordinate (G8) at (0.3045795679,0.1027824129);
        \coordinate (G9) at (0.3423101783,0.1348395019);
        \coordinate (G10) at (0.3766864961,0.1704699141);
        \coordinate (G11) at (0.4073716562,0.2093244979);
        \coordinate (G12) at (0.4340649725,0.2510225151);
        \coordinate (G13) at (0.4565048764,0.2951553655);
        \coordinate (G14) at (0.4744714788,0.3412905903);
        \coordinate (G15) at (0.4877887246,0.3889761095);
        \coordinate (G16) at (0.4963261177,0.4377446516);
        \coordinate (G17) at (0.5000000000,0.4871183324);
        \coordinate (G18) at (0.4987743711,0.5366133383);
        \coordinate (G19) at (0.4926612410,0.5857446666);
        \coordinate (G20) at (0.4817205123,0.6340308781);
        \coordinate (G21) at (0.4660593935,0.6809988154);
        \coordinate (G22) at (0.4458313482,0.7261882385);
        \coordinate (G23) at (0.4212345916,0.7691563353);
        \coordinate (G24) at (0.3925101477,0.8094820604);
        \coordinate (G25) at (0.3599394874,0.8467702619);
        \coordinate (G26) at (0.3238417674,0.8806555551);
        \coordinate (G27) at (0.2845706805,0.9108059180);
        \coordinate (G28) at (0.2425108300,0.9369260271);
        \coordinate (G29) at (0.1980728157,0.9587606064);
        \coordinate (G30) at (0.1516811082,0.9760992521);
        \coordinate (G31) at (0.1037132799,0.9887887995);
        \coordinate (G32) at (0.0541101801,0.9967658602);
        \coordinate (G33) at (0.0000000000,1.0000000000);
        \coordinate (G34) at (-0.0541101801,0.9967658602);
        \coordinate (G35) at (-0.1037132799,0.9887887995);
        \coordinate (G36) at (-0.1516811082,0.9760992521);
        \coordinate (G37) at (-0.1980728157,0.9587606064);
        \coordinate (G38) at (-0.2425108300,0.9369260271);
        \coordinate (G39) at (-0.2845706805,0.9108059180);
        \coordinate (G40) at (-0.3238417674,0.8806555551);
        \coordinate (G41) at (-0.3599394874,0.8467702619);
        \coordinate (G42) at (-0.3925101477,0.8094820604);
        \coordinate (G43) at (-0.4212345916,0.7691563353);
        \coordinate (G44) at (-0.4458313482,0.7261882385);
        \coordinate (G45) at (-0.4660593935,0.6809988154);
        \coordinate (G46) at (-0.4817205123,0.6340308781);
        \coordinate (G47) at (-0.4926612410,0.5857446666);
        \coordinate (G48) at (-0.4987743711,0.5366133383);
        \coordinate (G49) at (-0.5000000000,0.4871183324);
        \coordinate (G50) at (-0.4963261177,0.4377446516);
        \coordinate (G51) at (-0.4877887246,0.3889761095);
        \coordinate (G52) at (-0.4744714788,0.3412905903);
        \coordinate (G53) at (-0.4565048764,0.2951553655);
        \coordinate (G54) at (-0.4340649725,0.2510225151);
        \coordinate (G55) at (-0.4073716562,0.2093244979);
        \coordinate (G56) at (-0.3766864961,0.1704699141);
        \coordinate (G57) at (-0.3423101783,0.1348395019);
        \coordinate (G58) at (-0.3045795679,0.1027824129);
        \coordinate (G59) at (-0.2638644679,0.0746128255);
        \coordinate (G60) at (-0.2205644238,0.0506070612);
        \coordinate (G61) at (-0.1751078600,0.0310019184);
        \coordinate (G62) at (-0.1279693828,0.0159973810);
        \coordinate (G63) at (-0.0798158765,0.0057745326);
        \coordinate (G64) at (-0.0325860417,0.0005310661);
        \coordinate (GC) at (0,0.5);

        \path[fill=q64pale] (G1)

            -- (G2)
            -- (G3)
            -- (G4)
            -- (G5)
            -- (G6)
            -- (G7)
            -- (G8)
            -- (G9)
            -- (G10)
            -- (G11)
            -- (G12)
            -- (G13)
            -- (G14)
            -- (G15)
            -- (G16)
            -- (G17)
            -- (G18)
            -- (G19)
            -- (G20)
            -- (G21)
            -- (G22)
            -- (G23)
            -- (G24)
            -- (G25)
            -- (G26)
            -- (G27)
            -- (G28)
            -- (G29)
            -- (G30)
            -- (G31)
            -- (G32)
            -- (G33)
            -- (G34)
            -- (G35)
            -- (G36)
            -- (G37)
            -- (G38)
            -- (G39)
            -- (G40)
            -- (G41)
            -- (G42)
            -- (G43)
            -- (G44)
            -- (G45)
            -- (G46)
            -- (G47)
            -- (G48)
            -- (G49)
            -- (G50)
            -- (G51)
            -- (G52)
            -- (G53)
            -- (G54)
            -- (G55)
            -- (G56)
            -- (G57)
            -- (G58)
            -- (G59)
            -- (G60)
            -- (G61)
            -- (G62)
            -- (G63)
            -- (G64)
            -- cycle;


        \foreach \i/\j in {2/33,2/34,3/34,3/35,4/35,4/36,5/36,5/37,6/37,6/38,7/38}{\draw[diam] (G\i)--(G\j);}
        \foreach \i/\j in {7/39,8/39,8/40,9/40,9/41,10/41,10/42,11/42,11/43,12/43,12/44}{\draw[diam] (G\i)--(G\j);}
        \foreach \i/\j in {13/44,13/45,14/45,14/46,15/46,15/47,16/47,16/48,17/48,17/49,18/49}{\draw[diam] (G\i)--(G\j);}
        \foreach \i/\j in {18/50,19/50,19/51,20/51,20/52,21/52,21/53,22/53,22/54,23/54,23/55}{\draw[diam] (G\i)--(G\j);}
        \foreach \i/\j in {24/55,24/56,25/56,25/57,26/57,26/58,27/58,27/59,28/59,28/60,29/60}{\draw[diam] (G\i)--(G\j);}
        \foreach \i/\j in {29/61,30/61,30/62,31/62,31/63,32/63,32/64,33/64}{\draw[diam] (G\i)--(G\j);}

        \draw[borderedge] (G1)

            -- (G2)
            -- (G3)
            -- (G4)
            -- (G5)
            -- (G6)
            -- (G7)
            -- (G8)
            -- (G9)
            -- (G10)
            -- (G11)
            -- (G12)
            -- (G13)
            -- (G14)
            -- (G15)
            -- (G16)
            -- (G17)
            -- (G18)
            -- (G19)
            -- (G20)
            -- (G21)
            -- (G22)
            -- (G23)
            -- (G24)
            -- (G25)
            -- (G26)
            -- (G27)
            -- (G28)
            -- (G29)
            -- (G30)
            -- (G31)
            -- (G32)
            -- (G33)
            -- (G34)
            -- (G35)
            -- (G36)
            -- (G37)
            -- (G38)
            -- (G39)
            -- (G40)
            -- (G41)
            -- (G42)
            -- (G43)
            -- (G44)
            -- (G45)
            -- (G46)
            -- (G47)
            -- (G48)
            -- (G49)
            -- (G50)
            -- (G51)
            -- (G52)
            -- (G53)
            -- (G54)
            -- (G55)
            -- (G56)
            -- (G57)
            -- (G58)
            -- (G59)
            -- (G60)
            -- (G61)
            -- (G62)
            -- (G63)
            -- (G64)
            -- cycle;

        \draw[pendant] (G1)--(G33);
        \foreach \i in {2,...,32}{
            \node[dot] at (G\i) {};
        }
        \foreach \i in {34,...,64}{
            \node[dot] at (G\i) {};
        }
        \node[endpoint] at (G1) {};
        \node[endpoint] at (G33) {};
        \foreach \i in {2,...,32}{
            \node[font=\fontsize{7}{8}\selectfont, text=q64navy,
                  inner sep=0pt] at ($(GC)!1.17!(G\i)$) {\i};
        }
        \foreach \i in {34,...,64}{
            \node[font=\fontsize{7}{8}\selectfont, text=q64navy,
                  inner sep=0pt] at ($(GC)!1.17!(G\i)$) {\i};
        }
        \node[font=\fontsize{7}{8}\selectfont\bfseries,
              text=q64orange,inner sep=0pt]
              at ($(GC)!1.17!(G1)$) {1};
        \node[font=\fontsize{7}{8}\selectfont\bfseries,
              text=q64orange,inner sep=0pt]
              at ($(GC)!1.17!(G33)$) {33};

        \begin{scope}[shift={(0,-0.30)}]
            \draw[borderedge] (-0.54,0)--(-0.49,0);
            \node[font=\fontsize{7.5}{9}\selectfont,
                  text=q64navy,anchor=west] at (-0.478,0) {boundary};
            \draw[diam] (-0.16,0)--(-0.11,0);
            \node[font=\fontsize{7.5}{9}\selectfont,
                  text=q64navy,anchor=west] at (-0.098,0) {active diameters};
            \draw[pendant] (0.30,0)--(0.35,0);
            \node[font=\fontsize{7.5}{9}\selectfont,
                  text=q64navy,anchor=west] at (0.36,0) {pendant edge};
        \end{scope}
    \end{tikzpicture}
    \caption{The polygon $Q_{64}$ and its diameter graph.
Blue edges represent active unit distance constraints; the orange edge
$\{1,33\}$ is the pendant edge of its Graham type diameter graph.
The dark boundary follows the polygon's cyclic order, and the vertex labels
agree with the LSPectra input.}
    \label{fig:q64-diameter-framework}
\end{figure}

\begin{theorem}\label{thm:main}
	For each even \(n\) with \(6\le n\le\nmax\), let \(Q_n\) consist of the
	first \(2n\) components of the unique zero of \(\widehat F\) in the box
	\(Z_n\) recorded by LSPectra. Then:
    
	\begin{enumerate}[label=\textup{(\roman*)}, itemsep=0.6\baselineskip]
		\item
		\(Q_n\) is a strictly convex, positively oriented small \(n\)-gon.
		Its active diameter graph has edge set \(\mathcal E\) as defined in Section~\ref{sec:lspectra}. Figure~\ref{fig:diameter-graphs-all} shows the cases \(n=6,8,12,16,32\), and Figure~\ref{fig:q64-diameter-framework} shows \(Q_{64}\) with its pendant edge \(\{1,33\}\). Moreover,
		\[
		\delta(Q_n)\geq\delta_n,
		\qquad
		\Acal(Q_n)\in[A_n^-,A_n^+],
		\qquad
		A_n^+-A_n^-<10^{-11}.
		\]
		
		\item
		\emph{Linear independence:}
		\(Dg(Q_n)\) has full row rank, \(\rank Dg(Q_n)=|E(Q_n)|=n\), with
		\[
		s_{\min}(Dg(Q_n))\geq s_{n,L},
		\]
		and \(\Ical(Q_n)\cap R(Q_n)^\perp\) has dimension \(d=n-3\).
		
		\item
		\emph{Strict complementarity:}
		the multipliers of the multiplier identity are unique and satisfy
		\[
		\lambda_{n,L}\leq\lambda_{ij}\leq\lambda_{n,U}
		\qquad
		(\{i,j\}\in E(Q_n)).
		\]
		
		\item
		\emph{Reduced spectrum:}
		for the multipliers of \textup{(iii)}, set
		\[
		L_n=\Acal-
		\sum_{\{i,j\}\in E(Q_n)}\lambda_{ij}g_{ij}.
		\]
		The restriction of \(D^2L_n\) to
		\(\Ical(Q_n)\cap R(Q_n)^\perp\) has \(d\) simple eigenvalues
		\(\nu_1<\cdots<\nu_d\), enclosed in pairwise disjoint intervals,
		with \(\nu_d\le-\mu_n\). In particular,
		\[
		D^2L_n[v,v]\leq-\mu_n\lVert v\rVert^2
		\qquad
		\bigl(v\in\Ical(Q_n)\cap R(Q_n)^\perp\bigr).
		\]
		
		\item
		\emph{Quantitative stability:}
		every \(x\in(\R^2)^n\) satisfying
		\[
		\lVert x_i-x_j\rVert\leq1
		\quad (\{i,j\}\in E(Q_n)),
		\qquad
		\operatorname{dist}(x,Q_n)\leq\rho_n
		\]
		also satisfies
		\[
		\Acal(x)\leq
		\Acal(Q_n)-\kappa_n\operatorname{dist}(x,Q_n)^2.
		\]
		In particular, this holds for every small \(n\)-gon within
		distance \(\rho_n\) of \(Q_n\).
	\end{enumerate}
	
	In particular, \(Q_n\) is a strict local maximiser of \(\Acal\),
	modulo orientation preserving Euclidean motions, among all small
	\(n\)-gons, and admits no first order area increasing direction on
	its linearised feasible cone \(\Kcal(Q_n)\). The constants are
	those of Tables~\ref{tab:records} and~\ref{tab:cert-area}.
\end{theorem}

\begin{proof}[Proof of Theorem~\ref{thm:main}]
	Let \(n\) be even with \(6\le n\le\nmax\). By H1, \(\widehat F\) has a
	unique zero \((Q_n,\lambda,\xi)\) in \(Z_n\). Since H7 gives
	\(\Acal(Q_n)>0\), Lemma~\ref{lem:bordering-consistency} gives \(\xi=0\)
	and the multiplier identity for the pairs of \(\mathcal E\). As in
	Section~\ref{sec:lspectra}, \(g(Q_n)=0\) and H5 give
	\(\mathcal E=E(Q_n)\), so \(Q_n\) is small, and H8 makes it strictly
	convex and positively oriented. The interval evaluations over \(Z_n\)
	give the bounds in (i) and the multiplier bounds in (iii).
	
	The certification of \(s_{n,L}\) in Section~\ref{sec:lspectra} gives
	(ii). Since \(\rank Dg(Q_n)=n\), the space \(\ker Dg(Q_n)\) has
	dimension \(2n-n=n\) and contains \(R(Q_n)\), which has dimension \(3\)
	because \(\Acal(Q_n)>0\); hence the reduced space has dimension
	\(n-3\). Moreover \(Dg(Q_n)^\top\) is injective, so the multipliers
	are unique, which completes (iii).
	
	Since \(L_n\) is quadratic, its Hessian is constant and
	\(D^2L_n[v,v]\) is the constrained second variation of
	Theorem~\ref{thm:constrained-second-variation} at \(Q_n\). The verified inertia counts of
	Section~\ref{subsec:pipeline} increase from \(0\) to \(d\), by one
	across each recorded interval, and every reduced eigenvalue lies below
	the upper endpoint of the last interval, which is at most \(-\mu_n\);
	this encloses every reduced eigenvalue in its own interval, with
	multiplicity one, and gives (iv). Thus (ii)--(iv) allow us to apply
	Theorem~\ref{thm:quant-stability} with
	\[
	s_L=s_{n,L},\qquad
	\lambda_L=\lambda_{n,L},\qquad
	\lambda_U=\lambda_{n,U},\qquad
	\mu=\mu_n,
	\]
	and, by Remark~\ref{rem:exactness}, \(\Delta=3\).
	The exact rational evaluation of Section~\ref{sec:lspectra} gives the
	constants \(\rho_n\leq\rho\) and \(\kappa_n\leq\kappa\) of
	Table~\ref{tab:cert-area}, so the resulting estimate proves (v).
	
	Every small \(n\)-gon \(x\) with \(\lVert x-Q_n\rVert\le\rho_n\)
	satisfies \(\operatorname{dist}(x,Q_n)\le\rho_n\), so (v) gives
	\(\Acal(x)\leq\Acal(Q_n)\), with equality only on the Euclidean orbit
	of \(Q_n\). Finally, the positive multipliers in (iii)
	allow Proposition~\ref{prop:first-order-no-ascent} to exclude first order
	area increasing directions on \(\Kcal(Q_n)\).
\end{proof}

The records were produced with LSPectra~\cite{lspectra} version 0.24.5,
Julia~1.12.6 and IntervalArithmetic.jl~1.0.12. For every even \(n\) with
\(6\le n\le\nmax\), Krawczyk inclusion holds for the initial box \(Z\),
constructed with the same nominal radius \(10^{-13}\) in every component of 
\(u_0=(c_0,\lambda_0,\xi_0)\)
and H1--H8 hold with the bounds in Table~\ref{tab:records}.
Among these configurations, \(Q_{64}\) has the reduced eigenvalue closest
to zero; it remains more than \(9\tau\) below zero.

Three controls illustrate multiplicity, vanishing multipliers, and positive
second variation. The regular pentagon serves as a computational control for multiplicity and passes all eight verification checks H1--H8. Its reduced space has dimension \(2\), and the spectral verification encloses both reduced eigenvalues in a single interval with total multiplicity two.

For \(Q_4\), two multipliers vanish, so H4 fails, as
Proposition~\ref{prop:q4} predicts. Keeping all four diameters as equalities
gives reduced dimension \(d=1\). The spectrum \{\(-2,0,0\)\} in that proposition
instead concerns the larger equality space defined by the positive
multipliers.

For the regular hexagon, H1--H5, H7 and H8 hold. Its reduced spectrum has
two positive eigenvalues, counted with multiplicity, in an enclosure
containing \(\sqrt3/2\); the remaining reduced eigenvalues are negative.
The second order necessary condition for \(-\Acal\)
\cite[Thm.~12.5]{nocedalwright} therefore fails, and the hexagon is not a
local maximiser.

Table~\ref{tab:records} gives the bounds entering the stability estimate.
Here \(\delta_n\) bounds the inactive distance margin from below,
\(s_{n,L}\) bounds \(s_{\min}(Dg(Q_n))\) from below, and
\([\lambda_{n,L},\lambda_{n,U}]\) contains every active multiplier.
The quantity \(\mu_n\) is a lower bound for the negative of the largest
reduced eigenvalue of \(D^2L_n\) on \(\Ical(Q_n)\cap R(Q_n)^\perp\),
with its Euclidean inner product. Lower bounds are truncated downward
and upper bounds rounded upward at the last displayed digit.

For \(Q_4\), the zero eigenvalues in Proposition~\ref{prop:q4} reflect
the relaxed optimal family described there; the full diameter
inequalities restrict its feasible tangent directions to a proper cone.
For \(n=6,8,\ldots,\nmax\), the theorem establishes local maximality of
\(Q_n\) among all nearby small \(n\)-gons. The area of \(Q_6\) agrees to
six decimal places with Graham's value \(0.674981\ldots\)~\cite{graham},
but the theorem does not identify the configurations with global maximisers.

\subsection{Spectral Behaviour and Stability Bounds}

Table~\ref{tab:spectra} reports the smallest reduced eigenvalue
\(\nu_1\) and the two largest, \(\nu_{d-1}<\nu_d\), of
Theorem~\ref{thm:main}\textup{(iv)}. For each even \(n\), the full list of
intervals and their multiplicities is recorded under
\textit{evidence.reduced\_spectrum.eigenvalues} and the counts at their
endpoints under \textit{evidence.reduced\_spectrum.shifts} in the corresponding
LSPectra record~\cite{lspectra}. The eigenvalue enclosures give,
for even \(n\) with \(16\le n\le\nmax\),
\[
10.05<n|\nu_1|<10.68,
\qquad
1.13<\frac{n(\nu_d-\nu_{d-1})}{|\nu_d|}<2.14.
\]
The relative separation \((\nu_d-\nu_{d-1})/|\nu_d|\) decreases from
\(1.214\) at \(n=6\) to \(0.0178\) at \(n=\nmax\).

The spectrum also governs the stability constants of
Theorem~\ref{thm:main}\textup{(v)}: \(\kappa_n\) is obtained by rounding
\(\mu_n/8\) downward. Since \(\mu_n/(16\Delta)<\lambda_{n,L}\) for every
\(n\) in Table~\ref{tab:records}, the minimum in
Theorem~\ref{thm:quant-stability} is attained by the spectral term,
so that \(\rho_n\) is the downward rounding of \(s_{n,L}\mu_n/(48K)\),
with \(K\) evaluated at the bounds of Table~\ref{tab:records}.

The certified radius \(\rho_n\) is far below the lower bound
\(\delta_n\) for the inactive distance margin: the ratio
\(\rho_n/\delta_n\) decreases from about \(2.8\times10^{-3}\) at \(n=6\)
to below \(2\times10^{-6}\) at \(n=\nmax\).
Since \(K\ge q^2/\mu_n\), the formula for \(\rho_n\) gives
\[
\rho_n\le\frac{s_{n,L}\mu_n^2}{48q^2}.
\]
The small certified radii reflect the dependence of the
estimate on the reduced spectral gap. They are sufficient radii;
their sharpness is not established here.

\begin{remark}
    \label{rem:simple-competitors}
    By Remark~\ref{rem:shoelace}, \(\Acal\) is a signed area. Let
    \(D_n>0\) be the lower bound recorded by LSPectra for the determinants
    \(\det(z_{i+1}-z_i,z_j-z_i)\), \(j\notin\{i,i+1\}\), entering H8 over
    \(Z_n\). If \(\operatorname{dist}(x,Q_n)\le r\le\frac14\), choose a
    motion attaining this distance and let \(y\) be its image of \(x\).
    Then \(\lVert y_i-(Q_n)_i\rVert\le r\) for every \(i\); since every side
    and chord of \(Q_n\) has length at most \(1\), each of these
    determinants changes by at most \(4r+4r^2\le5r\) from \(Q_n\) to \(y\).
    For every \(n\), \(\rho_n<\frac14\) and the records give
    \(5\rho_n<D_n\); taking \(r=\rho_n\), \(y\) is
    strictly convex and positively oriented, and so is \(x\), since
    orientation preserving motions preserve both properties and the area.
    Thus every competitor in Theorem~\ref{thm:main}\textup{(v)} is strictly
    convex and positively oriented, and \(\Acal(x)\) is its enclosed area.
\end{remark}

\begin{table}[htbp]
	\caption{Bounds at the exact configurations \(Q_n\),
		\(n=6,8,\ldots,\nmax\), from the LSPectra~\cite{lspectra} records.}
	\label{tab:records}
	\centering
	\small
	\setlength{\tabcolsep}{6pt}
	\renewcommand{\arraystretch}{1.15}
	\begin{tabular}{lccccc}
		\toprule
		Polygon & \(\delta_n\) & \(s_{n,L}\) & \(\lambda_{n,L}\) & \(\lambda_{n,U}\) & \(\mu_n\)\\
		\midrule
		\(Q_{6}\) & \(0.2207794341\) & \(1.398\) & \(0.0415920318\) & \(0.1718857266\) & \(0.4628678899\)\\
		\(Q_{8}\) & \(0.1225787912\) & \(1.094\) & \(0.0380689324\) & \(0.1310708600\) & \(0.3064916518\)\\
		\(Q_{10}\) & \(0.0774263054\) & \(0.8883\) & \(0.0323841888\) & \(0.1055060193\) & \(0.2093529327\)\\
		\(Q_{12}\) & \(0.0532020485\) & \(0.7443\) & \(0.0277289987\) & \(0.0880806838\) & \(0.1436760844\)\\
		\(Q_{14}\) & \(0.0387602898\) & \(0.6394\) & \(0.0241082934\) & \(0.0755014536\) & \(0.09922242721\)\\
		\(Q_{16}\) & \(0.0294761687\) & \(0.5600\) & \(0.0212683681\) & \(0.0660206880\) & \(0.06986884663\)\\
		\(Q_{18}\) & \(0.0231616045\) & \(0.4979\) & \(0.0190000448\) & \(0.0586314357\) & \(0.05049135040\)\\
		\(Q_{20}\) & \(0.0186750376\) & \(0.4480\) & \(0.0171543828\) & \(0.0527163102\) & \(0.03745497467\)\\
		\(Q_{22}\) & \(0.0153743418\) & \(0.4072\) & \(0.0156270883\) & \(0.0478772927\) & \(0.02845825337\)\\
		\(Q_{24}\) & \(0.0128761290\) & \(0.3732\) & \(0.0143443148\) & \(0.0438469261\) & \(0.02208638069\)\\
		\(Q_{26}\) & \(0.0109401426\) & \(0.3444\) & \(0.0132528154\) & \(0.0404391427\) & \(0.01746316000\)\\
		\(Q_{28}\) & \(0.0094096665\) & \(0.3197\) & \(0.0123134460\) & \(0.0375206406\) & \(0.01403466754\)\\
		\(Q_{30}\) & \(0.0081789713\) & \(0.2983\) & \(0.0114968900\) & \(0.0349934931\) & \(0.01144208835\)\\
		\(Q_{32}\) & \(0.0071746436\) & \(0.2795\) & \(0.0107808091\) & \(0.0327841850\) & \(0.009447244997\)\\
		\(Q_{34}\) & \(0.0063444260\) & \(0.2630\) & \(0.0101479125\) & \(0.0308364748\) & \(0.007888337378\)\\
		\(Q_{36}\) & \(0.0056503058\) & \(0.2483\) & \(0.0095846235\) & \(0.0291066168\) & \(0.006653073725\)\\
		\(Q_{38}\) & \(0.0050640998\) & \(0.2352\) & \(0.0090801431\) & \(0.0275600811\) & \(0.005661980738\)\\
		\(Q_{40}\) & \(0.0045645559\) & \(0.2234\) & \(0.0086257799\) & \(0.0261692538\) & \(0.004857798607\)\\
		\(Q_{42}\) & \(0.0041354049\) & \(0.2127\) & \(0.0082144637\) & \(0.0249117912\) & \(0.004198595142\)\\
		\(Q_{44}\) & \(0.0037640221\) & \(0.2030\) & \(0.0078403868\) & \(0.0237694220\) & \(0.003653203273\)\\
		\(Q_{46}\) & \(0.0034404901\) & \(0.1941\) & \(0.0074987356\) & \(0.0227270613\) & \(0.003198140564\)\\
		\(Q_{48}\) & \(0.0031569311\) & \(0.1860\) & \(0.0071854881\) & \(0.0217721464\) & \(0.002815492821\)\\
		\(Q_{50}\) & \(0.0029070236\) & \(0.1785\) & \(0.0068972581\) & \(0.0208941311\) & \(0.002491436456\)\\
		\(Q_{52}\) & \(0.0026856477\) & \(0.1716\) & \(0.0066311750\) & \(0.0200840975\) & \(0.002215191340\)\\
		\(Q_{54}\) & \(0.0024886212\) & \(0.1652\) & \(0.0063847896\) & \(0.0193344537\) & \(0.001978268439\)\\
		\(Q_{56}\) & \(0.0023125008\) & \(0.1593\) & \(0.0061559997\) & \(0.0186386964\) & \(0.001773922275\)\\
		\(Q_{58}\) & \(0.0021544308\) & \(0.1538\) & \(0.0059429911\) & \(0.0179912228\) & \(0.001596747701\)\\
		\(Q_{60}\) & \(0.0020120270\) & \(0.1487\) & \(0.0057441901\) & \(0.0173871801\) & \(0.001442379648\)\\
		\(Q_{62}\) & \(0.0018832859\) & \(0.1438\) & \(0.0055582247\) & \(0.0168223442\) & \(0.001307267233\)\\
		\(Q_{64}\) & \(0.0017665144\) & \(0.1393\) & \(0.0053838937\) & \(0.0162930209\) & \(0.001188502161\)\\
		\bottomrule
	\end{tabular}
\end{table}

\begin{table}[htbp]
	\caption{Area enclosures \(\Acal(Q_n)\in[A_n^-,A_n^+]\) and
		stability constants for \(n=6,8,\ldots,\nmax\).
		The constants \(\rho_n\) and \(\kappa_n\) are rounded downward from
		the exact rational evaluation of Theorem~\ref{thm:quant-stability},
		using the bounds in Table~\ref{tab:records} and \(\Delta=3\)
		(Section~\ref{sec:lspectra}).}
	\label{tab:cert-area}
	\centering
	\small
	\setlength{\tabcolsep}{6pt}
	\renewcommand{\arraystretch}{1.15}
	\begin{tabular}{lcccc}
		\toprule
		Polygon & \(A_n^-\) & \(A_n^+\) & \(\rho_n\) & \(\kappa_n\)\\
		\midrule
		\(Q_{6}\) & \(0.674981442930\) & \(0.674981442931\) & \(6.152\times10^{-4}\) & \(0.05785848623\)\\
		\(Q_{8}\) & \(0.726868482751\) & \(0.726868482752\) & \(3.042\times10^{-4}\) & \(0.03831145647\)\\
		\(Q_{10}\) & \(0.749137345877\) & \(0.749137345878\) & \(1.506\times10^{-4}\) & \(0.02616911658\)\\
		\(Q_{12}\) & \(0.760729873448\) & \(0.760729873449\) & \(7.301\times10^{-5}\) & \(0.01795951055\)\\
		\(Q_{14}\) & \(0.767531011120\) & \(0.767531011121\) & \(3.516\times10^{-5}\) & \(0.01240280340\)\\
		\(Q_{16}\) & \(0.771861321980\) & \(0.771861321981\) & \(1.738\times10^{-5}\) & \(0.008733605828\)\\
		\(Q_{18}\) & \(0.774788165074\) & \(0.774788165075\) & \(8.976\times10^{-6}\) & \(0.006311418800\)\\
		\(Q_{20}\) & \(0.776858756001\) & \(0.776858756002\) & \(4.856\times10^{-6}\) & \(0.004681871833\)\\
		\(Q_{22}\) & \(0.778377330126\) & \(0.778377330127\) & \(2.746\times10^{-6}\) & \(0.003557281671\)\\
		\(Q_{24}\) & \(0.779524045160\) & \(0.779524045161\) & \(1.616\times10^{-6}\) & \(0.002760797586\)\\
		\(Q_{26}\) & \(0.780411119922\) & \(0.780411119923\) & \(9.860\times10^{-7}\) & \(0.002182895000\)\\
		\(Q_{28}\) & \(0.781111419017\) & \(0.781111419019\) & \(6.207\times10^{-7}\) & \(0.001754333442\)\\
		\(Q_{30}\) & \(0.781673926721\) & \(0.781673926722\) & \(4.019\times10^{-7}\) & \(0.001430261043\)\\
		\(Q_{32}\) & \(0.782132559738\) & \(0.782132559739\) & \(2.667\times10^{-7}\) & \(0.001180905624\)\\
		\(Q_{34}\) & \(0.782511405883\) & \(0.782511405884\) & \(1.811\times10^{-7}\) & \(0.0009860421722\)\\
		\(Q_{36}\) & \(0.782827954425\) & \(0.782827954426\) & \(1.254\times10^{-7}\) & \(0.0008316342156\)\\
		\(Q_{38}\) & \(0.783095152247\) & \(0.783095152249\) & \(8.851\times10^{-8}\) & \(0.0007077475922\)\\
		\(Q_{40}\) & \(0.783322749493\) & \(0.783322749494\) & \(6.349\times10^{-8}\) & \(0.0006072248258\)\\
		\(Q_{42}\) & \(0.783518202055\) & \(0.783518202056\) & \(4.622\times10^{-8}\) & \(0.0005248243927\)\\
		\(Q_{44}\) & \(0.783687290205\) & \(0.783687290206\) & \(3.412\times10^{-8}\) & \(0.0004566504091\)\\
		\(Q_{46}\) & \(0.783834551000\) & \(0.783834551002\) & \(2.550\times10^{-8}\) & \(0.0003997675705\)\\
		\(Q_{48}\) & \(0.783963585884\) & \(0.783963585885\) & \(1.928\times10^{-8}\) & \(0.0003519366026\)\\
		\(Q_{50}\) & \(0.784077283019\) & \(0.784077283020\) & \(1.474\times10^{-8}\) & \(0.0003114295570\)\\
		\(Q_{52}\) & \(0.784177980335\) & \(0.784177980337\) & \(1.138\times10^{-8}\) & \(0.0002768989175\)\\
		\(Q_{54}\) & \(0.784267586693\) & \(0.784267586695\) & \(8.866\times10^{-9}\) & \(0.0002472835548\)\\
		\(Q_{56}\) & \(0.784347673027\) & \(0.784347673029\) & \(6.969\times10^{-9}\) & \(0.0002217402843\)\\
		\(Q_{58}\) & \(0.784419541669\) & \(0.784419541671\) & \(5.522\times10^{-9}\) & \(0.0001995934626\)\\
		\(Q_{60}\) & \(0.784484279635\) & \(0.784484279637\) & \(4.409\times10^{-9}\) & \(0.0001802974560\)\\
		\(Q_{62}\) & \(0.784542799975\) & \(0.784542799977\) & \(3.542\times10^{-9}\) & \(0.0001634084041\)\\
		\(Q_{64}\) & \(0.784595874141\) & \(0.784595874144\) & \(2.866\times10^{-9}\) & \(0.0001485627701\)\\
		\bottomrule
	\end{tabular}
\end{table}

\begin{table}[htbp]
	\caption{Extreme reduced eigenvalues of \(D^2L_n\), \(n=6,8,\ldots,\nmax\), on
		\(\Ical(Q_n)\cap R(Q_n)^\perp\): the smallest \(\nu_1\) and the two largest
		\(\nu_{d-1}<\nu_d\), with \(\nu_d\le-\mu_n\). The three eigenvalues are rounded
		to eight decimals and lie within \(10^{-8}\) of their exact values. Their record enclosures have width at most \(2.6\times10^{-9}\). The last column gives
		\((\nu_d-\nu_{d-1})/|\nu_d|\), rounded to three decimals.}
	\label{tab:spectra}
	\centering
	\small
	\setlength{\tabcolsep}{6pt}
	\renewcommand{\arraystretch}{1.15}
	\begin{tabular}{lcccc}
		\toprule
		Polygon & \(\nu_1\) & \(\nu_{d-1}\) & \(\nu_d\) & \((\nu_d-\nu_{d-1})/|\nu_d|\)\\
		\midrule
		\(Q_6\) & \(-1.14208698\) & \(-1.02499156\) & \(-0.46286789\) & \(1.214\)\\
		\(Q_8\) & \(-1.09984120\) & \(-0.55866420\) & \(-0.30649165\) & \(0.823\)\\
		\(Q_{10}\) & \(-0.87475315\) & \(-0.30978807\) & \(-0.20935293\) & \(0.480\)\\
		\(Q_{12}\) & \(-0.81528786\) & \(-0.18464983\) & \(-0.14367609\) & \(0.285\)\\
		\(Q_{14}\) & \(-0.70226285\) & \(-0.11767682\) & \(-0.09922243\) & \(0.186\)\\
		\(Q_{16}\) & \(-0.62849940\) & \(-0.07920446\) & \(-0.06986885\) & \(0.134\)\\
		\(Q_{18}\) & \(-0.56585620\) & \(-0.05570488\) & \(-0.05049135\) & \(0.103\)\\
		\(Q_{20}\) & \(-0.51236906\) & \(-0.04060044\) & \(-0.03745498\) & \(0.084\)\\
		\(Q_{22}\) & \(-0.46970141\) & \(-0.03047353\) & \(-0.02845825\) & \(0.071\)\\
		\(Q_{24}\) & \(-0.43243247\) & \(-0.02344035\) & \(-0.02208638\) & \(0.061\)\\
		\(Q_{26}\) & \(-0.40095542\) & \(-0.01840835\) & \(-0.01746316\) & \(0.054\)\\
		\(Q_{28}\) & \(-0.37357115\) & \(-0.01471560\) & \(-0.01403467\) & \(0.049\)\\
		\(Q_{30}\) & \(-0.34965106\) & \(-0.01194572\) & \(-0.01144209\) & \(0.044\)\\
		\(Q_{32}\) & \(-0.32859420\) & \(-0.00982815\) & \(-0.00944725\) & \(0.040\)\\
		\(Q_{34}\) & \(-0.30989184\) & \(-0.00818200\) & \(-0.00788834\) & \(0.037\)\\
		\(Q_{36}\) & \(-0.29319377\) & \(-0.00688328\) & \(-0.00665307\) & \(0.035\)\\
		\(Q_{38}\) & \(-0.27818660\) & \(-0.00584509\) & \(-0.00566198\) & \(0.032\)\\
		\(Q_{40}\) & \(-0.26463156\) & \(-0.00500534\) & \(-0.00485780\) & \(0.030\)\\
		\(Q_{42}\) & \(-0.25232789\) & \(-0.00431886\) & \(-0.00419860\) & \(0.029\)\\
		\(Q_{44}\) & \(-0.24111100\) & \(-0.00375224\) & \(-0.00365320\) & \(0.027\)\\
		\(Q_{46}\) & \(-0.23084391\) & \(-0.00328046\) & \(-0.00319814\) & \(0.026\)\\
		\(Q_{48}\) & \(-0.22141133\) & \(-0.00288450\) & \(-0.00281549\) & \(0.025\)\\
		\(Q_{50}\) & \(-0.21271601\) & \(-0.00254973\) & \(-0.00249144\) & \(0.023\)\\
		\(Q_{52}\) & \(-0.20467508\) & \(-0.00226478\) & \(-0.00221519\) & \(0.022\)\\
		\(Q_{54}\) & \(-0.19721766\) & \(-0.00202072\) & \(-0.00197827\) & \(0.021\)\\
		\(Q_{56}\) & \(-0.19028265\) & \(-0.00181049\) & \(-0.00177392\) & \(0.021\)\\
		\(Q_{58}\) & \(-0.18381722\) & \(-0.00162841\) & \(-0.00159675\) & \(0.020\)\\
		\(Q_{60}\) & \(-0.17777537\) & \(-0.00146994\) & \(-0.00144238\) & \(0.019\)\\
		\(Q_{62}\) & \(-0.17211692\) & \(-0.00133137\) & \(-0.00130727\) & \(0.018\)\\
		\(Q_{64}\) & \(-0.16680660\) & \(-0.00120968\) & \(-0.00118850\) & \(0.018\)\\
		\bottomrule
	\end{tabular}
\end{table}

\FloatBarrier

\section*{Acknowledgements}
During the preparation of this article, we used OpenAI's
ChatGPT and Anthropic's Claude to explore mathematical arguments under the authors' guidance, with their questions and conjectures providing the starting point.
These exchanges helped us investigate and refine preliminary arguments and proof strategies, examine potential objections, and improve the exposition of the results. The assistance extended to the development and documentation of LSPectra and manuscript revision. It included the following activities:

\begin{itemize}

\item Proposing approaches to finding the optimal translation by
aligning the vertex centroids in the discussion preceding
Theorem~\ref{thm:quant-stability} and to deriving the explicit
constants \(\rho\) and \(\kappa\) in the theorem's quadratic area decrease estimate.

\item Examining whether the positivity assumption on the multipliers
could be removed from the local maximality criterion on the active
equality manifold (Proposition~\ref{prop:slice-local-max}).

\item Refining the computational framework, generating code under specific author instructions, and examining bottlenecks and possible improvements to performance and mathematical implementation. This included the input generator based on angular data from OPTIGON and verified approximate diagonalisation using Gershgorin bounds.

\item Generating software tests, implementing the Lean~4 formalisation of selected linear algebra results, and producing documentation from the existing code, including README files, comments, and docstrings.

\item Reviewing and improving the presentation of the quadrilateral analysis in Section~\ref{sec:models}, the computational procedures in Section~\ref{sec:lspectra}, and the results in Section~\ref{sec:main}. The quadrilateral analysis was written by the authors and subsequently revised with AI assistance.

\item Assisting with literature searches, figures,
\LaTeX{} preparation, and editorial revisions.

\end{itemize}

The authors take full responsibility for the mathematical arguments, computational results, and final manuscript, as well as for the accompanying software, tests, formalisation and documentation.
	

\end{document}